\input amstex
\documentstyle{amsppt}
\magnification=\magstep1

\pageheight{9.0truein}
\pagewidth{6.5truein}

\NoBlackBoxes
\TagsAsMath
\TagsOnLeft

\long\def\ignore#1{#1}

\ignore{
\input xy
\xyoption{matrix}\xyoption{arrow}\xyoption{curve}\xyoption{frame}
\def\edge{\ar@{-}}
}

\def\Rmod{R{}\operatorname{-Mod}}
\def\rmod{R{}\operatorname{-mod}}

\def\finpresR{{\operatorname{fin{.}pres-}}R}

\def\rad{\operatorname{rad}}

\def\dim{\operatorname{dim}}
\def\Ann{\operatorname{Ann}}
\def\Hom{\operatorname{Hom}}
\def\End{\operatorname{End}}
\def\Add{\operatorname{Add}}
\def\add{\operatorname{add}}

\def\supp{\operatorname{supp}}

\def\NN{{\Bbb N}}
\def\ZZ{{\Bbb Z}}
\def\RR{{\Bbb R}}
\def\A{{\Cal A}}
\def\B{{\Cal B}}
\def\C{{\Cal C}}

\def\K{{\Cal K}}
\def\CC{{\frak C}}
\def\DD{{\frak D}}

\def\AndFul{{\bf 1}}
\def\Aus{{\bf 2}}
\def\Auslarge{{\bf 3}}
\def\Ausmal{{\bf 4}}
\def\Azu{{\bf 5}}
\def\Azufac{{\bf 6}}
\def\BaeBruLen{{\bf 7}}
\def\BaeLen{{\bf 8}}
\def\Bal{{\bf 9}}
\def\Benson{{\bf 10}}
\def\BuHor{{\bf 11}}
\def\Chase{{\bf 12}}
\def\CohKap{{\bf 13}}
\def\Coh{{\bf 14}}
\def\CrawJon{{\bf 15}}
\def\CrawBoeI{{\bf 16}}
\def\CrawBoeII{{\bf 17}}
\def\EisGri{{\bf 18}}
\def\Fac{{\bf 19}}
\def\FacI{{\bf 20}}
\def\Fai{{\bf 21}}
\def\Fie{{\bf 22}}
\def\Fuchspaper{{\bf 23}}
\def\Fuchsbook{{\bf 24}}
\def\Ful{{\bf 25}}
\def\FuRei{{\bf 26}}
\def\Gara{{\bf 27}}
\def\Grif{{\bf 28}}
\def\Gru{{\bf 29}}
\def\GruJen{{\bf 30}}
\def\GruJenII{{\bf 31}}
\def\GruJenfinal{{\bf 32}}
\def\Her{{\bf 33}}
\def\Herz{{\bf 34}}
\def\ringswhose{{\bf 35}}
\def\Huistrong{{\bf 36}}
\def\HuiZimI{{\bf 37}}
\def\HuiZimexchange{{\bf 38}}
\def\HuiZimII{{\bf 39}}
\def\HuiZimIII{{\bf 40}}
\def\Hull{{\bf 41}}
\def\JenLenI{{\bf 42}}
\def\JenLen{{\bf 43}}
\def\JenZim{{\bf 44}}
\def\Kapbook{{\bf 45}}
\def\Kiel{{\bf 46}}
\def\KielSim{{\bf 47}}
\def\Koe{{\bf 48}}
\def\Kra{{\bf 49}}
\def\KraI{{\bf 50}}
\def\KraSao{{\bf 51}}
\def\Lara{{\bf 52}}
\def\Los{{\bf 53}}
\def\Mar{{\bf 54}}
\def\Myc{{\bf 55}}
\def\Okoh{{\bf 56}}
\def\Oso{{\bf 57}}
\def\Prest{{\bf 58}}
\def\PrestI{{\bf 59}}
\def\Rinspecies{{\bf 60}}
\def\Rin{{\bf 61}}
\def\RinTac{{\bf 62}}
\def\Schmid{{\bf 63}}
\def\Schulz{{\bf 64}}
\def\Sim{{\bf 65}}
\def\SimI{{\bf 66}}
\def\Simcox{{\bf 67}}
\def\SimArt{{\bf 68}}
\def\Simcounter{{\bf 69}}
\def\Simpad{{\bf 70}}
\def\Sten{{\bf 71}}
\def\Warf{{\bf 72}}
\def\Warfexchange{{\bf 73}}
\def\Warfnice{{\bf 74}}
\def\Weg{{\bf 75}}
\def\Zayed{{\bf 76}}
\def\Zel{{\bf 77}}
\def\Zimakad{{\bf 78}}
\def\Zimhab{{\bf 79}}
\def\Zim{{\bf 80}}
\def\Zimnote{{\bf 81}}
\def\ZimI{{\bf 82}}
\def\ZimII{{\bf 83}}
\def\Zwara{{\bf 84}}

\topmatter

\title  Purity, algebraic compactness, direct sum decompositions, and
representation type 
\endtitle

\rightheadtext{ Purity }

\author Birge Huisgen-Zimmermann \endauthor

\address Department of Mathematics, University of California, Santa Barbara, CA
93106, USA\endaddress
\email birge\@math.ucsb.edu\endemail

\thanks This work was partially supported by a grant from the National
Science Foundation (USA). \endthanks

\dedicatory Helmut Lenzing zu seinem 60\. Geburtstag gewidmet\enddedicatory

\endtopmatter

\document

\head  Introduction\endhead

The idea of purity pervades human culture, from
religion through sexual and moral codes, cleansing rituals and dietary
rules, to stratifications of societies into castes.  It does so to an
extent which is hardly backed up by the following definition found in the
Oxford English Dictionary: {\it freedom from admixture of any foreign substance
or matter\/}. 
In fact, in a process of
blending and fusing of ideas and ideals (which stands in
stark contrast to the dictionary meaning of the term itself), the concept 
has acquired dozens of additional connotations.  They have proved notoriously
hard to pin down, to judge by the incongruence of the emotional and
poetic reactions they have elicited.  Here is a small sample:    

\medskip  

{\it Purity is obscurity.\/} Ogden Nash.

{\it One cannot be precise and still be pure.\/} Marc Chagall.

{\it Unto the pure all things are pure.\/} New Testament.

{\it To the pure all things are impure.\/} Mark Twain.

{\it To the pure all things are indecent.\/} Oscar Wilde.

{\it Blessed are the pure in heart for they have so much more to talk about.\/}
Edith Wharton.

{\it Be thou as chaste as ice, as pure as snow, thou shalt not escape calumny.
Get thee to a nunnery, go.\/} W. Shakespeare.

{\it Necessary, for ever necessary, to burn out false shames and smelt the
heaviest ore of the body into purity.\/} D. H. Lawrence.

{\it Mathematics possesses not only truth, but supreme beauty -- a beauty cold
and austere, like that of a sculpture, without appeal to any part of our
weaker nature, sublimely pure, and capable of a stern perfection such as only
the greatest art can show.\/} Bertrand Russell. 

{\it Purity strikes me as the most mysterious of the virtues and the more I
think about it the less I know about it.\/} Flannery O'Connor.

\medskip

Even in the `sublimely pure' subject of mathematics, several notions of
purity were competing with each other a few decades ago.  However,
while the state of affairs is still rather muddled in contemporary society
at large, at least mathematicians have come to an agreement in the meantime. 
Section 1  traces the development of the concept within algebra in
rough strokes.  Due to the large number of contributions to the subject, we
can only highlight a selection.
	
The survey to follow should be read in conjunction with those of M.
Prest and G. Zwara (\cite{\PrestI} and \cite{\Zwara}). The first emphasizes the
model-theoretic aspect of the topic, while the second describes the impact of
generic modules on representation theory.  The generic modules form a
distinguished subclass of the class of pure injective (= algebraically compact)
modules, the importance of which became apparent through the seminal work of
Crawley-Boevey
\cite{\CrawBoeI}.  Since we wanted to minimize the overlap with an overview
article on  endofinite modules by the latter author, the article
\cite{\CrawBoeII} should be consulted as another supplement to our
survey.  For the sake of orientation, we offer the following diagram displaying
a hierarchy of noteworthy classes of algebraically compact modules.

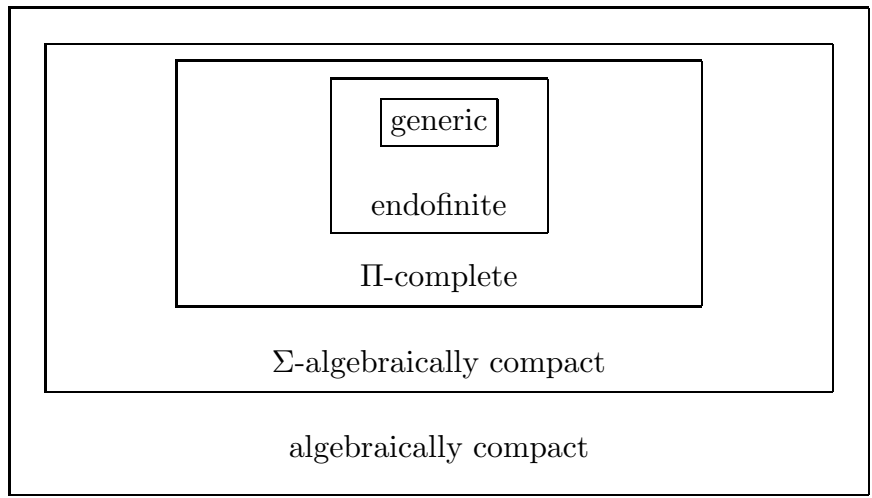
\ignore{
$$\xymatrixcolsep{2pc}\xymatrixrowsep{1pc}
\xymatrix{
\save[0,0]+(-4,0);[4,4]+(4,-3)**\frm{-}\restore 
\save[0,0]+(-8,4);[5,4]+(8,-5)**\frm{-}\restore &&&&\\
 & \save[0,0]+(0,6);[2,2]+(0,-3)**\frm{-}\restore &\txt{generic}
\save**\frm{-}\restore \save[0,0]+(-12,4);[1,0]+(12,-4)**\frm{-}\restore &&\\
 &&\txt{endofinite} &&\\
 &&\Pi\hbox{-complete} &&\\
 &&\Sigma\hbox{-al\-ge\-bra\-i\-cal\-ly compact} &&\\
 &&\hbox{algebraically compact} &&  
}$$
}  

In the present overview, we  mainly address the outer layers of the
stack. We  conclude with a  discussion of product-complete modules,
the study of which  provides a natural bridge to the class of
endofinite modules.  The interest of the latter class also becomes clear as one
approaches it from alternate angles, e.g., via the duality theory
developed in
\cite{\CrawBoeII}. Zwara's survey picks up the story at this point and zeroes in
on the generic modules, i.e., the non-finitely generated indecomposable
endofinite modules.  The extensive spread of applications of these concepts is
further illustrated by the contributions of D. Benson and H. Krause to this
volume (\cite{\Benson} and \cite{\KraI}):  Benson encounters them in his study of
phantom maps between modules over group algebras, while Krause uses them towards
new characterizations of tame representation type.  

What follows below is almost exclusively expository. 
However, a few results rounding off the picture (e.g., the characterization of
purity in terms of matrix groups in Proposition 4) appear to be new, as are
several of our arguments.  Moreover, in some instances, we streamline results
which can essentially be found in the literature by taking the present state of
the art into account.

The red thread that will lead us through this discussion is provided in
Section 2, in the form of a number of global decomposition problems going
back to work of Koethe (1935) and Cohen-Kaplansky (1951).  Long before
the impact of the subject on the representation theory of tame finite
dimensional algebras  surfaced, these problems had put a spotlight on
($\Sigma$-)al\-ge\-bra\-ic compactness and linked it to finite
representation type.  Section 3  contains the most important
characterizations of the algebraically compact and the
$\Sigma$-al\-ge\-bra\-i\-cal\-ly compact modules, as well as the functorial underpinnings
on which they are based.  Section 4  takes us back to the decomposition
problems stated at the outset, and Section 5  evaluates the outcome in
representation-theoretic terms. The topic of product-completeness, addressed in
Section 6, supplements both the decomposition theory of direct powers of modules
as described in Theorem 10 of Section 4, and the impact of `globally nice'
decompositions of products on the representation type of the underlying ring
discussed in Section 5. The very short final section should be seen as an
appendix, meant to trigger further research in a direction that has been somewhat
neglected in the recent past.

\head 1. Purity and algebraic compactness -- definitions and a brief
history\endhead

It was already in the first half of this century that the concept of a pure
subgroup of an abelian group proved pivotal in accessing the structure,
first of $p$-groups, then also of torsionfree and mixed abelian groups, as well
as of modules over PID's.  Given a PID $R$, a submodule
$A$ of an $R$-module $B$ is called {\it pure\/} if, for all $r \in R$, the
intersection $A\cap rB$ equals $rA$.  As is well-known, the submodule $T(B)$
consisting of the torsion elements of  $B$ is always pure, and
in exploring the $p$-primary components of $T(B)$, major headway is
gained by studying a tell-tale class of pure submodules of the simplest
possible structure, namely the `basic' ones.  They are determined up to
isomorphism by the following requirements:  Given a prime $p \in R$, a
submodule $B_0$ of a $p$-primary $R$-module $B$ is called {\it basic\/} if it is
a direct sum of cyclic groups which is pure in $B$ and has the additional
property of making
$B/B_0$ divisible.  In fact, a natural extension of this concept to arbitrary
modules over discrete valuation domains provides us with one of the reference
points in the classification of the algebraically compact abelian groups
sketched below. 
  
In the 1950's, it became apparent that suitable variations of the original
notion of purity should yield important generalizations of split embeddings
in far more general contexts, and several extensions of the concept, naturally
all somewhat akin in spirit, appeared in the literature. Finally, in 1959,
Cohn's definition \cite{\Coh} won the `contest'.  Given left modules $A$ and $B$
over an arbitrary associative ring $R$, Cohn called a monomorphism
$f:A\rightarrow B$ {\it pure\/} in case tensoring with any 
right $R$-module $X$ preserves injectivity in the induced map $id \otimes f:
X\otimes_R A \rightarrow X \otimes_R B$;  by extension, a  short
exact sequence $0 \rightarrow A \rightarrow B \rightarrow C \rightarrow 0$ is
labeled {\it pure exact\/} in case the monomorphism from $A$ to $B$ is pure.
When, in 1960, Maranda \cite{\Mar} undertook a study of the modules $M$ for
which the contravariant
$\Hom$-functor
$\Hom_R(-,M)$ takes maps from certain restricted classes of monomorphisms in
$\Rmod$ to epimorphisms in the category of abelian groups, the
class which quickly became the most popular in this game was that of pure
monomorphisms. Accordingly, modules $M$ such that
$\Hom_R(-,M)$ preserves exactness in pure exact sequences were labelled {\it
pure injective\/}. 

One of the fundamental demands on a good concept of purity is this: It should
yield convenient criteria for recognizing direct summands  --  the philosophy
being that purity is a first step toward splitting.  One hopes
for readily verifiable conditions which take a pure inclusion the rest of the
way to a split one.  Here are two sample statements of this flavor for Cohn's
purity;   we will re-encounter them in Section 3.
  
\proclaim{Observations 0}

$\bullet$  If $R$ is a Dedekind domain and $M$ a pure submodule of an
$R$-module $N$ such that $\text{Ann}_R(M) \ne 0$, then $M$ is a direct
summand of $N$. 

$\bullet$ If $R$ is a left perfect ring and $Q$ a pure submodule of
a projective left $R$-module $P$, then $Q$ is a direct summand of $P$.
\endproclaim

\demo{Proof}  The first statement will arise as an immediate consequence of
Theorem 6 in Section 3.  To establish the second we will show that,
given a pure exact sequence
$$0 \rightarrow Q \rightarrow P \rightarrow P/Q \rightarrow 0$$
 of left modules over a ring $R$, such that $P$
is  flat, the quotient $P/Q$ is flat as
well.  This will yield our claim since a perfect base ring will make flat
modules projective.

To check flatness of $P/Q$, let $0 \rightarrow A \rightarrow B \rightarrow C
\rightarrow 0$ be any short exact sequence of right $R$-modules, and
consider the following diagram which has exact rows and columns due to our
setup.

\ignore{
$$\xymatrixcolsep{4pc}\xymatrixrowsep{2pc}
\xymatrix{
 &&0 \ar[d] \ar@{-->}[r] &{\hbox{Ker}}(\phi) \ar[d] 
\ar@{-->}@'{ @+{[0,0]+(20,0)} @+{[1,0]+(20,5)} @+{[1,0]+(20,-6)}
@+{[1,-2]+(-20,-5)} @+{[2,-2]+(-20,-2)} @+{[3,-2]+(-20,0)} }
[3,-2]  \\  
0 \ar[r] &A\otimes_R Q \ar[r] \ar[d] &A\otimes_R P \ar[r] \ar[d]
&A\otimes_R (P/Q) \ar[r] \ar[d]^{\phi} &0\\
0 \ar[r] &B\otimes_R Q \ar[r] \ar[d] &B\otimes_R P \ar[r] \ar[d] &B\otimes_R
(P/Q) \ar[r] \ar[d] &0\\
 &C\otimes_R Q \ar[r]^{\psi} \ar[d] &C\otimes_R P \ar[r] \ar[d] &C\otimes_R
(P/Q) \ar[r] \ar[d] &0\\
 &0 &0 &0
}$$
}

\noindent Applying the Snake Lemma and using the fact that the map $\psi$ is
injective by hypothesis, we deduce that $\text{Ker}(\phi) = 0$ as required. \qed
\enddemo

The theory of purity and pure injectivity could also be baptized {\it the
theory of systems of linear equations for modules\/}. A typical
such system has the following format: Starting with a left $R$-module $M$, an
element
$(m_i)_{i\in I}\in M^I$, and a row-finite matrix $(r_{ij})_{i\in I,j\in
J}$ of elements of $R$, one considers the system
$$ \qquad \qquad \sum_{j\in J} r_{ij}X_j =m_i \qquad\qquad (i\in I), \tag
\dagger$$
 and calls each element in $M^J$ which satisfies all equations of
$(\dagger)$ a solution in
$M$.  A first indication of a connection between such systems and the theory of
purity and pure injectivity was exibited by Fieldhouse (see \cite{\Fie}), who
observed that purity of a submodule $A$ of an $R$-module $B$ is equivalent to
the following equational condition: {\sl Every finite system of linear equations
with right-hand sides in $A$, which is solvable in $B$, has a solution in
$A$\/}.  This result makes one anticipate a bridge between linear systems and
pure injectivity as well.  Such a bridge does in fact exist, but relies on
systems that need not be finite.  Here it is:

\proclaim{Theorem 1} {\rm (Warfield, 1969 \cite{\Warf})} For $M\in\Rmod$, the
following statements are equivalent:

{\rm (1)} $M$ is pure injective.

{\rm (2)} Any system of the form $(\dagger)$ which is finitely solvable in $M$
(i.e\., has the property that, for any finite subset $I' \subseteq I$, the
finite subsystem of equations labeled by $i \in I'$ is solvable in $M$) has a
global solution in $M$. 

{\rm (3)} $M$ is a direct summand of a compact Hausdorff $R$-module $N$ (the
latter is to mean that $N$ is a compact Hausdorff abelian group such that all
multiplications by elements of $R$ are continuous). \endproclaim

In this
instance, Warfield acted primarily as a coordinator of results and ideas from
various parts of the literature, most of the implications having been
previously known in a variety of specialized contexts, some in full
generality.  It is quite enlightening to follow the line of successive
modifications of the pertinent concepts.  From work of Kaplansky
\cite{\Kapbook}, \L o\'s \cite{\Los}, and Balcerzyk \cite{\Bal} done in the
fifties, through papers of Mycielski
\cite{\Myc}, Butler-Horrocks
\cite{\BuHor}, Kie\l pi\'nski \cite{\Kiel}, Weglorz \cite{\Weg}, Fuchs
\cite{\Fuchspaper}, and Stenstr\o m \cite{\Sten} scattered over the sixties,
it incrementally leads  to the present theory.  

Warfield's coordination was very
much called for, as the obvious interest of the topological condition
$(3)$ of Theorem 1 had triggered various exploratory trips in its own right. 
In 1954, Kaplansky published a characterization of those abelian groups which
arise as algebraic direct summands of compact Hausdorff groups, referring to
them as `algebraically compact' \cite{\Kapbook}. Originally, he had set out
to describe the compact abelian groups, but realized that, from an algebraic
viewpoint, the former class allowed for a cleaner description and appeared 
 more natural. In particular, he provided a complete classification of the
algebraically compact abelian groups: They are precisely the 
direct sums of divisible groups and groups which are Hausdorff and complete in
their
$\ZZ$-adic topologies. The latter can be written uniquely as direct products of
factors $A_p$ which are  complete and Hausdorff in their $p$-adic topologies,
respectively, with $p$ tracing the primes.  Each $A_p$, in turn, can be
canonically viewed as a module over the ring of $p$-adic integers and pinned
down up to isomorphism in terms of its basic submodule, a direct sum of cyclic
p-goups and copies of the p-adic integers, thus leading to a convenient full
set of invariants. Subsequently, equivalent descriptions of the class of groups
exhibited by Kaplansky were given by Mycielski, \L o\'s, Weglorz, and Fuchs (in
the setting of noetherian rings).  Most notable were the descriptions in terms
of `equational compactness' conditions. The first general concept of
algebraic compactness  was introduced by
Mycielski in 1964, for arbitrary algebraic systems in fact. Restricted to
modules, it just amounts to condition $(2)$ of Theorem 1.  For the sake of
emphasis, we repeat: 

\definition{Definition} Given an associative ring $R$ with identity, a left
$R$-module $M$ is called {\it algebraically compact\/} if any system
$$\sum_{j\in J} r_{ij}X_j =m_i \qquad\qquad (i\in I),$$ based on a row-finite
matrix
$(r_{ij})$ with entries in $R$ and $m_i\in M$, which is finitely
solvable in $M$, has a global solution in $M$. 

Moreover, $M$ is said to be {\it $\Sigma$-al\-ge\-bra\-i\-cal\-ly compact\/} in
case all direct sums of copies of $M$ are algebraically compact.
\enddefinition 

So, in particular, Theorem 1 tells us that the pure injective and the
algebraically compact modules coincide.  We will give the latter terminology
preference in the sequel. 

\definition{\it Remark\/} It is of course natural to also wonder about
characterizations of pure projective modules, i.e., of those modules which are
projective relative to pure exact sequences. Such characterizations are much
more readily obtained than useful descriptions of their pure injective
counterparts. According to Warfield
\cite{\Warf} and Fieldhouse \cite{\Fie}, the pure projective modules are
precisely the direct summands of direct sums of finitely presented modules.
\enddefinition

With little effort, one obtains the following first list of examples of
algebraically compact (alias pure injective) modules that goes beyond the most
obvious ones, the injectives; a second installment of examples can be found in
part B of Section 3.
\smallskip

1. \cite{\Warf} Any module $M$ can be purely embedded into a pure injective
module, namely its `Bohr compactification', as follows:
$$M\rightarrow \Hom_\ZZ(\Hom_\ZZ(M,\RR/\ZZ),\RR/\ZZ),$$ the assignment being
evaluation. If one equips $\Hom_\ZZ(M,\RR/\ZZ)$ with the discrete topology and
the  `double dual' of $M$ with the compact-open topology, one thus arrives
at a compact Hausdorff $R$-module. As a consequence, one can construct pure
injective resolutions of a module $M$ and measure the deviation of
$M$ from pure injectivity by means of its `pure injective dimension'.

2. \cite{\Warf} If $R$ is a commutative local noetherian domain with
maximal ideal
$\frak m$, complete in its $\frak m$-adic topology, then $R$ is
algebraically compact (as an $R$-module). In particular, this is true for the
ring of $p$-adic integers.  Note, however, that the ring of $p$-adic integers
fails to be $\Sigma$-al\-ge\-bra\-i\-cal\-ly compact (both as an abelian group and as a
module over itself).

3. \cite{\Fuchsbook} If $R$ is a countable ring, $A$ any left $R$-module, and
$\Cal F$ a non-principal ultrafilter on $\NN$, then the ultrapower
$A^{\NN}/{\Cal F}$ is algebraically compact.

4. Every artinian module over a commutative ring is $\Sigma$-al\-ge\-bra\-i\-cal\-ly
compact (see the examples following Theorem 6 below).

5. If $R$ is a commutative artinian principal ideal ring, then all $R$-modules
are direct sums of cyclic submodules.  The latter being finite in number, up to
isomorphism, Example 4 shows all objects of $\Rmod$ to be algebraically
compact in that case.
\smallskip

We will sketch a proof for the fifth remark (apparently folklore), since it
is part of the red thread that will lead us through this survey. By the
Chinese Remainder Theorem, $R$ is a finite direct product of local rings, and
hence it is harmless to assume that $R$ is a local artinian principal ideal
ring with maximal ideal $\frak m$ say.  Now all ideals of $R$ are powers
of $\frak m$, and each homomorphism from a power $\frak m^n$ to $R$ sends
$\frak m^n$ back to $\frak m^n$ and can therefore be extended to $R$  --  just
use the fact that $\frak m$ is principal. This shows that $R$ is self-injective;
in fact, $R$ being an artinian ring, $R$ is
$\Sigma$-injective as an $R$-module.  Now let $M$ be any nonzero module.  We
provide a decomposition of the desired ilk by induction on the least
natural number
$N$ such that $\frak m^{N+1} M
= 0$.  Clearly $M$ is a module over the ring $R / \frak m ^{N+1}$ which, as
we just saw, is $\Sigma$-injective over itself.  Let $(x_i)_{i\in I}$ be a
maximal family of elements of $M$ such that the sum of the $Rx_i$ is direct,
with each of the $Rx_i$ isomorphic to $R$.  Due
to injectivity, the sum  $\sum_{i \in I} Rx_i$ is then a direct summand of $M$,
say $M = \sum_{i \in I} Rx_i \oplus M'$, and due to the maximal choice of our
family, $M'$ does not contain a copy of $R / \frak m ^{N+1}$, i.e\.,
$M'$ is annihilated by $\frak m^N$.  Our claim follows by induction.
\medskip

Point 5 shows, in particular, that all artinian principal ideal rings have
finite representation type in the sense to follow.  In fact, it has long been
known that, among the commutative artinian rings, the  principal ideal
rings are precisely the ones having this property.

\definition{Definition}  A ring $R$ is said to have {\it finite representation
type\/} if it is left artinian and if, up to isomorphism, there are only
finitely many indecomposable finitely generated left $R$-modules.
\enddefinition

Finite representation type is actually left-right symmetric,
as was shown by Eisenbud and Griffith in \cite{\EisGri}.

\head 2. A global decomposition problem \endhead

For better focus, we interject two
problems which, on the face of it, are  only loosely  connected
with our main theme.  The connection turns out to be much closer than
anticipated at first sight.  In fact, these problems have motivated a major
portion of the subsequent work on algebraic compactness. 

\proclaim{Global Problems} {\rm (Koethe \cite{\Koe}, Cohen-Kaplansky
\cite{\CohKap})} For which rings
$R$ is every right $R$-module

{\rm (a)} a direct sum of finitely generated modules?

{\rm (b)} a direct sum of indecomposable modules?
\endproclaim

Work on the commutative case was initiated by Koethe in 1935,
continued by Cohen-Kaplansky in 1951, and completed by Griffith in 1970
for part (a), and by Warfield in 1972 for part (b). 

\proclaim{Theorem 2} {\rm (\cite{\Koe, \CohKap, \Grif, \Warfnice})}
For any
commutative ring $R$,  conditions {\rm (a)} and {\rm (b)} above are
equivalent and satisfied if and only if $R$ is an artinian principal ideal
ring. \endproclaim

As we noted at the end of Section 1, among the commutative
rings, the artinian principal ideal rings are precisely the ones having
finite representation type.  Moreover, we observed that these rings enjoy the
property that all their modules are algebraically compact. As we will see in
Theorem 13 of Section 4, this property in turn characterizes the artinian
principal ideal rings, which rounds off the `commutative solution' to our
problems.  The reasons for the most interesting
implications, namely that either of the two conditions (a), (b) forces
$R$ to be an artinian principal ideal ring, can be roughly summarized as
follows: In general, large direct products of modules exhibit a high
resistance to infinite direct sum decompositions; more precisely, in most cases,
infinite direct sum decompositions of direct products can be traced back to
infinite decompositions of finite sub-products.

In the noncommutative situation, the problems become far more challenging. In
his 1972 paper, Warfield stated:  ``For non-commutative rings, the questions
raised in this paper seem to be much more difficult.  All that seems to be known
is that any ring satisying [the above conditions (a), (b)] is necessarily
[left] artinian."  At present, there is still at least one link missing to a
truly satisfactory resolution of these problems. The main key to what is known
is an equivalent characterization of
$\Sigma$-al\-ge\-bra\-ic compactness, which will be presented in the next
section.

\head 3. Characterizations of ($\Sigma$-)algebraically compact modules\endhead

We begin with a subsection dedicated to the main technical resource of the
subject, introduced independently by Gruson-Jensen \cite{\GruJen} and
W. Zimmermann \cite{\Zim}.  It is the more general functorial framework of
\cite{\Zim} and
\cite{\HuiZimI} which we will describe below, since the extra generality
adds to the transparency of the arguments.

\subhead{A\. Product-compatible functors and matrix functors}\endsubhead

\definition{Definition} {\rm (1)} A $p$-{\it functor\/} on $\Rmod$ is a
subfunctor $P$ of the forgetful functor $\Rmod \rightarrow \bold{Ab}$ which
commutes with direct products, i.e\., $P$ assigns to each left $R$-module $M$
a subgroup
$PM$ such that $f(PM) \subseteq PN$ for any homomorphism $f: M \rightarrow N$,
and
$P \bigl(\prod_{i \in I}M_i \bigr) = \prod_{i \in I}(PM_i)$ for any direct
product
$\prod_{i \in I}M_i$ in $\Rmod$.

Note that any $p$-functor automatically commutes with direct sums (this being
actually true for any subfunctor of the forgetful functor).

{\rm (2)} A {\it pointed matrix\/} over $R$ is a row-finite matrix $\A =
(a_{ij})_{i \in I, j \in J}$ of elements from $R$, paired with a column index
$\alpha \in J$.  Given a pointed matrix $(\A,\alpha)$, we call the following
$p$-functor $[\A,\alpha]$ on $\Rmod$ a {\it matrix-functor\/}:  For any
$R$-module $M$, the subgroup $[\A, \alpha]M$ is defined to be the $\alpha$-th
projection of the solution set in $M$ of the homogeneous system
$$\sum_{j \in J} a_{ij} X_j = 0 \qquad \qquad \text{for all}\ \ i \in I;$$
in other words, 
$$[\A, \alpha]M = \{m \in M \mid \exists\ \text{a solution}\ (m_j) \in M^J
\text{\ of the above system with}\ m_{\alpha} = m\}.$$
Further, we call $[\A,\alpha]$ a {\it finite matrix functor\/} in case the
matrix $\A$ is finite.
\enddefinition

Given a (finite) matrix functor $[\A,\alpha]$ on $\Rmod$ and a left $R$-module
$M$, we will call the subgroup $[\A,\alpha]M$  a {\it (finite) matrix
subgroup\/} of $M$.  The finite matrix subgroups were labeled ``sousgroupes
de d\'efinition fini'' by Gruson and Jensen, and ``$pp$-definable subgroups'' by
the model theorists, Prest, Herzog, Rothmaler, and others.

\proclaim{Observation 3. Basic properties of matrix functors} The first two
give alternate descriptions of matrix functors.

{\rm (1)} $[\A, \alpha] = \Hom_R(Z, -)(z)$ for a suitable left 
$R$-module
$Z$ and $z \in Z$; conversely, every functor of the form $\Hom_R(Z, -)(z)$ is a
matrix functor for a suitable matrix $\A$.  The finite matrix functors are
precisely the functors $\Hom_R(Z, -)(z)$ with finitely presented $Z$.

{\rm (2)} The finite matrix subgroups of $M$ moreover coincide with
the kernels of the $\ZZ$-linear maps $M \rightarrow Z \bigotimes_R M$,
$m \mapsto z\otimes m$, where $Z$ is a finitely presented right $R$-module
and $z \in Z$.  

{\rm (3)} The class of matrix functors is closed under arbitrary intersections
and finite sums.  The finite matrix functors are closed under
\underbar{finite} intersections and finite sums.

{\rm (4)} If $M$ is an $R$-$S$-bimodule, then every matrix subgroup of $_RM$ is
an $S$-submodule of $M$.\endproclaim

In the context of part (4), the most important candidates for $S$ will be
the opposite of the endomorphism ring of $M$, as well as subrings of the center
of $R$.  How easily matrix functors can be manipulated
is evidenced by the easy proofs of the above observations; we include one
sample argument to make our point.

\demo{Proof of part of {\rm (3)}}  We will show that finite sums of matrix
functors are again of that ilk.  So let $(\A = (a_{ij})_{i \in I, j \in J},
\alpha)$ and
$(\B = (b_{kl})_{k \in K, l \in L}, \beta)$ be two pointed matrices with
entries in $R$.  It is clearly harmless to assume that the occurring index sets
are all disjoint.  Let $(\C = (c_{uv})_{u \in U, v\in V}, \gamma)$ be
defined as follows:  Assuming that $\gamma$ belongs to none of the sets
$I,J,K,L$, we set $U = \{\gamma\} \cup I \cup K$ and $V = \{\gamma\} \cup J \cup
L$, and define $c_{uv}$ via $c_{\gamma \gamma} = 1$, $c_{\gamma \alpha} =
c_{\gamma \beta} = -1$, $c_{ij} = a_{ij}$ whenever $(i,j) \in I \times J$,
$c_{kl} = b_{kl}$ whenever $(k,l) \in K \times L$, and $c_{uv} = 0$ in all
other cases.  It is straightforward to check that, for any left $R$-module
$M$, we have $[\A,\alpha]M + [\B,\beta]M = [\C,\gamma]M$.  Moreover, we remark
that $\C$ is finite if $\A$ and $\B$ are. \qed \enddemo

\proclaim{Important instances of matrix functors}
 
 $\bullet$ Whenever $\frak a$ is a finitely generated right ideal of
$R$, the assignment $M \mapsto \frak a M$ defines a finite matrix functor on
$\Rmod$.  

$\bullet$  For every subset $T$ of $R$, the assignment $M \mapsto
\text{Ann}_M(T)$ is a matrix functor; it is finite in case $T$ is.

More generally: Whenever $\frak a$ is a finitely generated right ideal
and $T$ a subset of $R$, the conductor $(\frak a M:T)$ is a matrix subgroup of
$M$.

$\bullet$ All finitely generated $\text{End}_R(M)$-submodules of
$M$ are matrix subgroups. \endproclaim

Pinning down matrices for the first three types of examples is
straightforward.  To verify the last:  By Observation 3(1), each cyclic
$\text{End}_R(M)$-submodule of $M$ is a matrix subgroup; now use
Observation 3(3) to move to finite sums.

As a first application of matrix functors to our present objects of interest,
we will give a characterization of purity in terms of finite matrix
subgroups.  Ironically, this description provides a perfect formal
parallel to that of purity over PID's, while all of the initial
attempts to generalize the concept of purity, which were based on formal
analogy (with requirements such as $\frak a N \cap M = \frak a M$ for all
cyclic or finitely generated right ideals of the base ring), were eventually
discarded as not quite strong enough to best serve their purpose.

\proclaim{Proposition 4} A submodule $M$ of a left $R$-module $N$ is pure if
and only if $[\A,\alpha]N \cap M = [\A,\alpha]M$ for all finite matrix
functors $[\A,\alpha]$ on $\Rmod$.  
\endproclaim

\demo{Proof} The proof for pure left exactness of finite matrix functors is
straightforward (cf\. \cite{\Zimhab}).

Now assume that $M \subseteq N$ satisfies the above intersection property.
To verify purity of the inclusion, let 
$$ \qquad \sum_{j=1}^s a_{ij} X_j = m_i  \qquad \qquad  (1 \le i \le
r) \tag \dagger$$
 be a finite linear system with $m_i \in M$,  which is solvable in
$N$.  We prove its solvability in $M$ by induction on $r$.  For $r=1$, let $\A =
(-1,a_{11},\dots, a_{1s})$ be a single row with the first column labeled
$\alpha$; then $m_1 \in [\A,\alpha]N \cap M = [\A,\alpha]M$ by construction and
hypothesis, which provides a solution of $(\dagger)$ in $M$.  Next suppose
that
$r\ge 2$, pick a solution
$(v_j)$ of $(\dagger)$ in $N$, and use the induction hypothesis to procure a
solution $(u_j)$ in $M$ of the first $r-1$ equations.  Define $m_r' : =
\sum_{j=1}^s a_{rj}(v_j - u_j) =  m_r- \sum_{j = 1}^s a_{rj} u_j \in M$, and
observe that $\sum_{j=1}^s a_{ij}(v_j - u_j) = 0$ for $i = 1,
\dots, r-1$.  We infer that $m_r'$ belongs to $[\B,\beta]N \cap M$, where
$\B$ is the matrix
$$\pmatrix 0 & a_{11} & a_{12} & \cdots & a_{1s}\\  
\vdots & \vdots & \vdots & &
\vdots \\ 0 & a_{r-1,1} & a_{r-1,2} &
\cdots & a_{r-1,s} \\ -1 & a_{r1} & a_{r2} & \cdots & a_{rs} \endpmatrix,$$
and $\beta$ labels the first column of $\B$.  By hypothesis, we thus
have $m_r' \in [\B,\beta]M$, which provides us with a family $(u_j')$
in $M$ satisfying $\sum_{j=1}^s a_{ij} u_j' = 0$ for $1 \le i \le r-1$
and $\sum_{j=1}^s a_{rj} u_j' =  m_r'$.  This yields a solution $(u_j +
u_j')$ of $(\dagger)$ in $M$ as required. \qed 
\enddemo

Moreover, numerous properties of a ring $R$ linked to the behavior of direct
products in $\Rmod$ can be conveniently recast in terms of matrix
functors.  We give two examples taken from \cite{\Zimhab, \Zim}:

 $\bullet$ A left $R$-module $M$ is flat if and only if $[\A, \alpha]M
= [\A, \alpha]R \cdot M$ for arbitrary finite matrix functors $[\A, \alpha]$. 
One deduces that all direct sums of copies of a flat module $M$ are flat if
and only if, for each finite matrix functor
$[\A,
\alpha]$, there exists a finitely generated right ideal $\frak a \subseteq
[\A, \alpha]R$ with $[\A, \alpha]M = \frak a M$.

 $\bullet$ Thus: flatness is inherited by arbitrary direct products of
flat left $R$-modules (i.e., $R$ is right coherent) if and only if all
finite matrix subgroups of $_RR$ are finitely generated right ideals.

\subhead{B. $p$-functors and algebraic compactness}
\endsubhead

The following equivalent descriptions of
algebraically compact and $\Sigma$-al\-ge\-bra\-i\-cal\-ly compact modules will 
prove extremely useful in the sequel. In each of the two theorems, the
equivalence of
$(1)$ and $(3)$ was independently established by Gruson-Jensen and
Zimmermann (\cite{\GruJenII, \Zim}); condition $(2)$  --  a significant
strengthening of the necessary side, in light of the third of the above
instances of matrix subgroups  --  was added by the latter. 

\proclaim{Theorem 5}  The following statements are equivalent for any left
$R$-module $M$:

{\rm (1)} $M$ is algebraically compact.

{\rm (2)} Every family of residue classes
 $(m_l+P_lM)_{l\in L}$ with $m_l\in M$ and $p$-functors $P_l$, which
has the finite intersection property, has non-empty intersection.

{\rm (3)} Same as {\rm (2)} with $P_l$ replaced by finite matrix functors. 
\endproclaim

In the special situation of a Pr\"ufer domain $R$ (i.e., a commutative
semihereditary integral domain), a specialized version of this result had
already been obtained by Warfield \cite{\Warf}:  In that case, the matrix
subgroups of the form
$(\frak a M:T)$, where $\frak a$ is a finitely generated ideal and $T$ a
finite subset of $R$, are representative.

Since Zimmermann's elementary proof of this theorem \cite{\Zim, Satz
2.1} does not exist in English translation, we  include a sketch below.
We begin with a remark which will come in handy in other contexts as
well.  Namely, locally  --  i.e., on the closure of a given module $M$ under
direct sums and products in $\Rmod$  --  any
$p$-functor acts like a matrix functor.  Indeed, given a $p$-functor $P$ on
$\Rmod$, we write
$PM = \{m_i \mid i \in I \}$, let $M_i$ be a copy of $M$ for each $i$, and set 
$\overline m = (m_i)_{i \in I} \in \prod_{i \in I} M_i$.  Then $PM =
\Hom_R(\prod_{i \in I}M_i, M)(\overline m)$, and Observation 3(1) yields our
claim.  

\demo{Proof of Theorem 5}  `$(1) \implies (2)$': Assuming (1), we start with
a family of residue classes $(m_l + P_l M)_{l \in L}$ having the
finite intersection property.  By the preceding remark, we may assume that
the $P_l$ are matrix functors, say $P_l = [\A_l, \alpha_l]$ with $\A_l = (a_{ij}^l)_{i
\in I_l, j \in J_l}$.  We construct
a system of equations which is finitely solvable in
$M$, and a global solution to which will provide us with an element in the
intersection of our family. It is clearly harmless to assume that all of the
index sets $I_l$ and $J_l$ are pairwise disjoint. Choose an element $\alpha$
contained in none of these.  We set
$I =
\bigcup_{l \in L} I_l$, $J = \{\alpha\}\cup
\bigcup_{l \in L} J_l \setminus \{\alpha_l\}$, and define a row-finite matrix
$\A = (a_{ij})_{i \in I, j \in J}$ as follows:  $a_{i\alpha} = a^l_{i \alpha_l}$
if $i \in I_l$;
$a_{ij} = a_{ij}^l$ if $i \in I_l$ and $j \in J_l \setminus \{\alpha_l\}$;
finally we set $a_{ij} = 0$ in all other cases.  It is immediate that
$[\A, \alpha]M = \bigcap_{l \in L}[\A_l, \alpha_l]M$.  The matrix $\A$ will
serve as coefficient matrix of our system.  As for its right-hand side, we
define
$m \in M^I$ by stringing up the elements $m_i = a_{i \alpha} m_l$ for $i \in
I_l$ in the only plausible fashion.  It is then straightforward to check that
the system
$\sum_{j \in J} a_{ij}X_j = m_i$ for
$i \in I$ is finitely
solvable by construction.  Hence it has a global solution, say $(z_i)_{i
\in I}$, and one readily verifies that $z_{\alpha}$ belongs to the
intersection of the family of residue classes with which we started out.

`$(3) \implies (1)$':  This time, we begin with a system 
$$ \qquad \sum_{j \in J} a_{ij}X_j = m_i \qquad \text{for} \qquad i
\in I, \tag \dagger$$
 which is finitely solvable in
$M$.  The crucial step of our argument is the following easy
consequence of (3) and finite solvability of $(\dagger)$: Namely, for each index
$\alpha
\in J$, there exists an element $y_{\alpha} \in M$ such the the system 
$$\sum_{j \in J \setminus \{\alpha\}} a_{ij}X_j = m_i - a_{i \alpha}
y_{\alpha}$$
is again finitely solvable, as follows. Given any finite subset $I' \subseteq
I$, let $y_\alpha(I')$ be the $\alpha$-th component of a solution of the finite
system
$\sum_{j\in J} a_{ij}X_j= m_i$ for $i\in I'$, and let $\A(I')$ be the matrix
consisting of the rows of $\A$ labelled by $I'$. Then the family
$\bigl( y_\alpha(I') + [\A(I'),\alpha]M \bigr)$, where $I'$ runs through the
finite subsets of $I$, has the finite intersection property. Therefore its
intersection is nonempty by (3), and any element $y_\alpha$ in this
intersection satisfies our requirement.

 Next we consider the set
$\K$ of all pairs
$(K, y)$, where
$K$ is a subset of $J$ and $y = (y_k) \in M^K$ is such that the system 
$$ \qquad  \sum_{j \in J \setminus K} a_{ij}X_j \ \ = \ \ m_i -
\sum_{k \in K} a_{ik} y_k \qquad \text{for} \qquad i \in I $$
is in turn finitely solvable in $M$.  We equip this set of pairs with the
standard order, namely $(K,y) \le (K',y')$ if $K
\subseteq K'$ and
$y_k = y_k'$ whenever $k \in K$.  Clearly, $\K \ne \varnothing$.  One checks
that the set $\K$ is inductively ordered, and denotes by $(K_0,z)$ a maximal
element of $\K$.  Our initial statement, applied to the finitely solvable system
$\sum_{j
\in J\setminus K_0} a_{ij}X_j = m_i - \sum_{k \in K_0}a_{ik} z_{k}$,
now yields $J = K_0$, which makes $z$ a global solution of $(\dagger)$.
\qed
\enddemo

\proclaim{Theorem 6}  The following statements are equivalent for any left
$R$-module $M$:

{\rm (1)} $M$ is $\Sigma$-al\-ge\-bra\-i\-cal\-ly compact.

{\rm (2)} Every countable descending chain
 $P_1\supseteq P_2\supseteq P_3\supseteq
\cdots$ of $p$-functors becomes stationary  
on $M$.

{\rm (3)} Same as {\rm (2)} with $P_l$ replaced by finite matrix functors.
\endproclaim

We point out that conditions $(2)$ and $(3)$ could just as well have been
phrased as follows:  `$M$ has the descending chain condition for
$p$-functorial subgroups' (or, equivalently, `$M$ has the descending chain
condition for finite matrix subgroups'), since the classes of $p$-functors and
finite matrix functors are both closed under finite intersections.  Indeed,
given a descending chain $P_1 M \supseteq P_2 M \supseteq \dots$ of
$p$-functorial subgroups of $M$, we obtain $P_i M = Q_i M$ for the descending
chain $Q_i = P_1 \cap \dots \cap P_i$ of $p$-functors. 

\demo{Proof of Theorem 6}  In view of Theorem 5, the implication `$(3)
\implies (1)$' is just an analogue of the much older result that artinian
modules are linearly compact
\cite{\Zel}; it is left as an exercise.  In the arguments given by
Zimmermann and Gruson-Jensen, the novel implication `$(1) \implies (2)$' is
obtained via a detour through direct products, which we sketch because it
pinpoints the importance of product-compatibility of the functors we are
considering:  Clearly,
$\Sigma$-al\-ge\-bra\-ic compactness of $M$ forces the natural (pure) embedding of
the direct sum $M^{(\NN)}$ in the direct product $M^{\NN}$ to split, say
$\prod_{n \in \NN} M_n = C \oplus \bigoplus_{n \in \NN} M_n$, where each
$M_n$ is a copy of $M$.  Let
$\pi: \prod_{n \in \NN} M_n \rightarrow C$ be the corresponding projection,
and $\pi_n: \prod_{i \in \NN} M_i \rightarrow M_n$ the canonical maps. 
Moreover, assume we have a descending chain  $P_1 \supseteq P_2 
\supseteq P_3
\supseteq \dots$ of
$p$-functors with $m_n \in P_n M_n \setminus P_{n+1} M_n$ for
$n \in \NN$.  Set $x = (m_n)$, $y_n = (m_1, \dots, m_n, 0, \dots)$, $x_n =
x - y_n$, and decompose $x$ and $x_n$ in the form $x = s + c$ and $x_n = s_n
+ c_n$ with $s, s_n \in \bigoplus_{n \in \NN} M_n$ and $ c, c_n \in C$;
then, clearly $s = y_n + s_n$.  Since $x_n$ belongs to $\prod_{n \in
\NN} P_{n+1} M_n =  P_{n+1} \bigl( \prod_{n \in \NN}M_n \bigr)$ by
construction,
$s_n$ belongs to $P_{n+1} \bigl(\bigoplus_{n \in \NN}M_n \bigr)$.  So if we
choose an index
$N \in \NN$ with
$\pi_N(s) = 0$, we infer that $m_N = \pi_N(y_N) = \pi_N (s - s_N) = -
\pi_N(s_N)$ lies in
$P_{N+1}M_N$, which contradicts our choice of $m_N$. \qed \enddemo

In Section 4, these concepts and results will turn out tailored to measure for
tackling the problems advertized in Section 2.  However, the two preceding
theorems lend themselves to a large number of further applications.  We
present a small selection, starting with a few additional examples; most of
them are obvious in light of the above theory, the others we tag with
references.

\definition{Further examples of ($\Sigma$-)algebraically compact modules}

$\bullet$  Whenever $M$ is an $R$-$S$-bimodule which is artinian
over $S$, $M$ is $\Sigma$-al\-ge\-bra\-i\-cal\-ly compact as an $R$-module. 
This is, for instance, true if $M$ is an $R$-module which is artinian over its
endomorphism ring or over the center of $R$.  In particular, for any
Artin algebra $R$, the category $\rmod$ consists entirely of
$\Sigma$-al\-ge\-bra\-i\-cal\-ly compact modules.

 $\bullet$ Suppose that $M$ is a
$(\Sigma)$-al\-ge\-bra\-i\-cal\-ly compact left $R$-module, $P$ a $p$-functor on
$\Rmod$ and $S$ a subring of $R$ such that $PM$ is an $S$-submodule of $M$. 
Then both $PM$ and $M/PM$ are ($\Sigma$-)al\-ge\-bra\-i\-cal\-ly compact over
$S$.  In particular:  If $\frak a$ is an ideal of $R$ which is  finitely
generated on the right, then $\frak a M$ and $M/
\frak a M$ are ($\Sigma$-)al\-ge\-bra\-i\-cal\-ly compact (see
\cite{\HuiZimI}). 

 $\bullet$  If $M$ is a module over a Dedekind domain $R$, then
$M$ is $\Sigma$-al\-ge\-bra\-i\-cal\-ly compact if and only if $M = M_1 \oplus
M_2$ where $M_1$ is divisible and $M_2$ has nonzero $R$-annihilator. 

 $\bullet$  If $M$ is finitely generated over a commutative
noetherian ring, then $M$ is $(\Sigma)$-al\-ge\-bra\-i\-cal\-ly compact
precisely when $M$ is artinian. (For the nontrivial implication, see
\cite{\Zim}.)   

 $\bullet$ Suppose $S$ is a left $\Sigma$-al\-ge\-bra\-i\-cal\-ly
compact ring (e.g\. a field), $(X_i)_{i \in I}$ a family of independent
indeterminates over
$S$, and $k$ any positive integer.  Then the truncated polynomial ring 
$$R = S[X_i \mid i \in i]/(X_i \mid i \in i)^k$$ is in turn
$\Sigma$-al\-ge\-bra\-i\-cal\-ly compact as a left module over itself (see
\cite{\HuiZimI}).   

$\bullet$ If $S$ is left algebraically compact, then so is every power
series ring $R = S[[X_i \mid i \in I]]$ ([loc\. cit.]).  Note, however, that $R$
fails to be
$\Sigma$-al\-ge\-bra\-i\-cal\-ly compact in case $I \ne \varnothing$.

 $\bullet$ If $S$ is any ring and $G$ a group, then the group ring
$SG$ is left algebraically compact if and only if $S$ has this property and $G$
is finite (see \cite{\Zimnote}).

 $\bullet$ Suppose $M$ is an $R$-$S$-bimodule.  If $A$ is an
algebraically compact left
$R$-module, then the left $S$-module
$\Hom_R(M,A)$ is algebraically compact as well -- see below.  (This
`compactification by passage to a suitable dual' is akin to the
Pontrjagin dual.)  In particular, the endomorphism ring $S$ of any
algebraically compact left module is algebraically compact on the right, the
change of side being due to the convention that $S$ acts on the left of
$M$.
\enddefinition

Since the last of the above remarks, concerning the passage of algebraic
compactness from
$_R A$ to $_S{}\Hom_R(M, A)$, will be relevant in the sequel, we include a
justification based on the adjointness of Hom and tensor product:  Start with a
pure monomorphism $U \rightarrow V$ of left
$S$-modules, and observe that the pure injectivity of $A$ forces the upper row
(and hence also the lower one) of the following commutative diagram to be an
epimorphism:

\ignore{
$$\xymatrixcolsep{4pc}\xymatrixrowsep{2pc}
\xymatrix{
{\ \hbox{Hom}_R(M\otimes_S V,\, A)\ } \ar[r] \ar[d]^{\cong}
&{\ \hbox{Hom}_R(M\otimes_S U,\, A)\ } \ar[d]^{\cong}\\
\hbox{Hom}_S(V,\, \hbox{Hom}_R(M,A)) \ar[r] &\hbox{Hom}_S(U,\,
\hbox{Hom}_R(M,A))
}$$
}

The following equivalence is due to Faith \cite{\Fai}.

\proclaim{Corollary 7}   For any injective left $R$-module $M$, the following
statements are equivalent:

{\rm (1)} $M$ is $\Sigma$-injective (meaning that every direct sum of copies
of $M$ is injective).

{\rm (2)}  R has the ascending chain condition for $M$-annihilators,
i.e., for left ideals of the form
$\text{Ann}_R(U)$, where
$U$ is a subset of $M$.
\endproclaim

\demo{Proof} Note that condition (1) is tantamount to $\Sigma$-al\-ge\-bra\-ic
compactness of $M$ under our hypotheses, since $M^{(I)}$ is pure in $M^I$ for
any set $I$. The matrix subgroups of an injective module
$M$ are precisely the annihilators in
$M$ of subsets of
$R$ (prove $[\A, \alpha]M = \Ann_M \Ann_R [\A, \alpha]M$ or consult \cite{\Zim,
Beispiel 1.1}).  Thus Theorem 6 shows
$(1)$ to be equivalent to the descending chain condition for
$R$-annihilators in
$M$.  But the latter descending chain condition has an equivalent flip side,
namely the ascending chain conditon for
$M$-annihilators in $R$. \qed \enddemo

\proclaim{Corollary 8}  Any pure submodule of a $\Sigma$-al\-ge\-bra\-i\-cal\-ly
compact module is a direct summand. 
\endproclaim

\demo{Proof}  Let $M$ be $\Sigma$-al\-ge\-bra\-i\-cal\-ly compact and $U\subseteq M$ a
pure submodule.  By Proposition 4, $[\A, \alpha]M \cap U = [\A, \alpha]U$ for
every finite matrix functor $[\A, \alpha]$, and hence $U$ inherits the
descending chain condition for finite matrix subgroups from $M$.  Therefore,
$U$ is in turn $\Sigma$-al\-ge\-bra\-i\-cal\-ly compact and consequently
a direct summand of $M$. \qed
\enddemo 

The next result is part of the bridge taking us back to our global
decomposition problems.  We will obtain the first part as another consequence
of Theorem 6, while the second can be most elegantly derived from the
following category equivalence due to Gruson-Jensen
\cite{\GruJenfinal} which provides a very general tool for extending
properties of injective modules to algebraically compact ones:  Let
$\finpresR$ be the full subcategory of
$\operatorname{mod}$-$R$ consisting of the finitely presented right $R$-modules,
and let 
$(\finpresR,\bold {Ab})$ be the category of all additive covariant
functors from
$\finpresR$ to the abelian goups.  This functor category is a
Grothendieck category, and the functor 
$$\Rmod \rightarrow (\finpresR,\bold {Ab})$$ 
defined by $M\mapsto
- \otimes_R M$ induces a category equivalence from the full subcategory of 
algebraically compact objects of
$\Rmod$ to the full subcategory of injective objects of  
$(\finpresR,\bold {Ab})$. 

\proclaim{Proposition 9} {\rm (Further assets of
($\Sigma$-)al\-ge\-bra\-i\-cal\-ly compact modules)}

{\rm (1)} Every
$\Sigma$-al\-ge\-bra\-i\-cal\-ly compact module is a direct sum of indecomposable
summands with local endomorphism rings.
 
{\rm (2)} If $M$ is an algebraically compact left module with endomorphism ring
$S$, then $S/\rad(S)$ is von Neumann regular and right self-injective; moreover,
idempotents can be lifted modulo
$\rad(S)$.  

In particular, $S$ is local in case $M$ is indecomposable.

{\rm (3)}  Each strongly invariant submodule $M$ of an algebraically compact
module has the exchange property, i.e., can be shifted inside direct sum
grids as follows: Given any equality of left $R$-modules $M \oplus X =
\bigoplus_{i\in I} Y_i$, there exist submodules $Y_i' \subseteq Y_i$ such that
$\bigoplus_{i \in I} Y_i = M \oplus \bigoplus_{i \in I} Y_i'$.

(Here a submodule $M$ of a module $N$ is said to be strongly
invariant in case each homomorphism $\phi$ from $M$ to $N$ satisfies
$\phi(M) \subseteq M$.)       
\endproclaim

\demo{Proof}
Part (1) \cite{\HuiZimI}: Suppose that $M$ is
$\Sigma$-al\-ge\-bra\-i\-cal\-ly compact, and let $(U_i)_{i \in I}$ be a maximal family
of independent indecomposable submodules of $M$ such that $U = \bigoplus_{i
\in I} U_i$ is pure in $M$.  By Corollary 8, $U$ is a direct summand of $M$, say
$M= U \oplus V$.  If $V$ were nonzero, we could pick a nonzero element $v \in
V$ and choose a pure submodule $W \subseteq V$ which is maximal with the
property of excluding $v$.  Again we would obtain splitting, $V = W \oplus
Z$, which would provide us with an indecomposable summand $Z$ of $V$.  But the
existence of such a $Z$ is incompatible with the maximal choice of the family
$(U_i)$. That the
$U_i$ even have local endomorphisms rings, will follow from part (2).  

For part (2), we use the above category equivalence and the
well-known fact that endomorphism rings of injectives have the listed
properties; indeed, the argument given by Osofsky in
\cite{\Oso} for modules carries over to Grothendieck categories.  (Note
however that the direct argument given in \cite{\HuiZimI} is a bit slicker
than the original proof for injective modules, one of the reasons being
that algebraic compactness is passed on to endomorphism rings by the last
of the examples following Theorem 6.)

As for part (3), we will only sketch a proof for the finite exchange
property and refer the reader to \cite{\HuiZimexchange} for a complete
argument.  Let
$M$ be a strongly invariant submodule of an algebraically compact left
$R$-module
$N$, and let
$S$ be the endomorphism ring of $M$.  From the last of the remarks following
Theorem 6, we know that
$S = \Hom_R(M,N)$ is an algebraically compact right $S$-module.  Therefore
$S$ is an exchange ring in the sense of Warfield by part (2), which
easily implies the finite exchange property of $M$ (see
\cite{\Warfexchange}).      
\qed
\enddemo

To conclude the section, we point to two examples to mark what we consider
potentially dangerous curves:  Namely, while any right artinian ring
is necessarily $\Sigma$-al\-ge\-bra\-i\-cal\-ly compact on the left, on the right it
need not even be algebraically compact.  An example demonstrating this was
given by Zimmermann in \cite{\Zimnote}.  Moreover, while the descending
chain condition for finite matrix subgroups of a module $M$ entails the
descending chain condition for arbitrary matrix subgroups, the analogous
implication fails for the ascending chain condition as the following example
shows:  One of the consequences we derived from Theorem 6 says that
the  algebra $R = K[X_i \mid i \in \NN]/(X_i \mid i \in
\NN)^2$ over a field $K$ is $\Sigma$-al\-ge\-bra\-i\-cal\-ly compact.  Hence, if $E$
denotes the minimal injective cogenerator for $R$, then
 the $R$-module $E = \Hom_R(R,E)$ has the ascending chain condition for finite
matrix subgroups (use Tool 17 below to see this).  Note that the ascending chain
condition for {\it arbitrary} matrix subgroups would amount to $E$ being
noetherian over its endomorphism ring $S$, since all finitely generated
$S$-submodules of $E$ are in fact matrix subgroups.  But the $S$-module $E$
fails to be noetherian:  Indeed, if it were, then the equality $A = \Ann_R
\Ann_E (A)$ for any ideal $A$ of $R$ would make $R$ artinian, an obvious
absurdity.  On the other hand:

\definition{\it Remark\/}  Suppose that $M$ is a direct sum of finitely
presented objects in $\rmod$.  Then the $R$-module $M$ satisfies the ascending
chain condition for finite matrix subgroups if and only if it is noetherian
over its endomorphism ring.

To see this, let $M = \bigoplus_{i \in I}M_i$ with finitely presented summands
$M_i$, and let $S$ be the endomorphism ring of $M$.   
For the nontrivial implication, suppose that $M$ has the ascending chain condition
for finite matrix subgroups.   Clearly, $M$ is noetherian over $S$ in case all
finitely generated $S$-submodules are finite matrix subgroups.  By Observation
3(3), it thus suffices to prove that each cyclic $S$-submodule $Sm$
of $M$ is a finite matrix subgroup.  Choose a finite subset $I' \subseteq I$
such that $m \in \bigoplus_{i \in I'} M_i$.  Then $Sm =
\Hom_R\bigl(\bigoplus_{i \in I'} M_i,M \bigr)(m)$ is indeed a finite matrix
subgroup of $M$ by Observation 3(1). \qed
\enddefinition

To learn about the module-theoretic impact of the maximum condition
for finite matrix subgroups, consult \cite{\ZimI}, for generalizations of
classical results to modules satisfying this condition, see \cite{\ZimII}.

\head{4. Return to our global decomposition problems} \endhead

As we mentioned in Section 2, the arguments settling the commutative case rest
on the fact that certain large direct products tend to resist nontrivial direct
sum decompositions.  It is therefore not surprising that, also in the
non-commutative situation, the decomposition properties of such
products are crucial.  

\proclaim{Theorem 10}  For a left $R$-module $M$, the following conditions are
equivalent:

{\rm (1)} There exists a cardinal number $\aleph$ such that every direct
product of copies of $M$ is a direct sum of $\aleph$-generated modules.

{\rm (2)} Every direct product of copies of $M$ is a direct sum of submodules
with local endomorphism rings.

{\rm (3)} $M$ is $\Sigma$-al\-ge\-bra\-i\-cal\-ly compact.
\endproclaim

The implication `$(1) \implies (3)$' is due to Gruson-Jensen
\cite{\GruJenII}, `$(2) \implies (3)$' to the author of this article
\cite{\ringswhose}, whereas `$(3) \implies (1)$' can already be found in work
of Kie\l pi\'nski going back to 1967 \cite{\Kiel}.  The implication `$(3)
\implies (2)$' is a consequence of Proposition 9.  Credit should also go to
Chase, since the ideas he developed in \cite{\Chase} play an essential role in
the proofs.  These ideas can be extracted and upgraded to the following format:

\proclaim{Lemma 11} {\rm (Chase)}  Let
$$f: \prod_{i\in\NN} U_i \longrightarrow \bigoplus_{j\in J} V_j$$ be a
homomorphism of $R$-modules, and $P_1\supseteq P_2\supseteq P_3\supseteq
\cdots$ a descending chain of $p$-functors. Then there exists a natural number
$n_0$ such that
$$f\biggl( P_{n_0} \prod_{i\ge n_0} U_i \biggr) \qquad \subseteq \qquad
\bigoplus_{\text{finite}} V_j \ \ +\ \ \bigcap_{n\in\NN} P_n\bigl( \bigoplus_J
V_j \bigr).$$
\noindent (Here $\bigoplus_{\text{finite}}$ stands for a direct sum extending
over some finite subset of $J$.)
\endproclaim

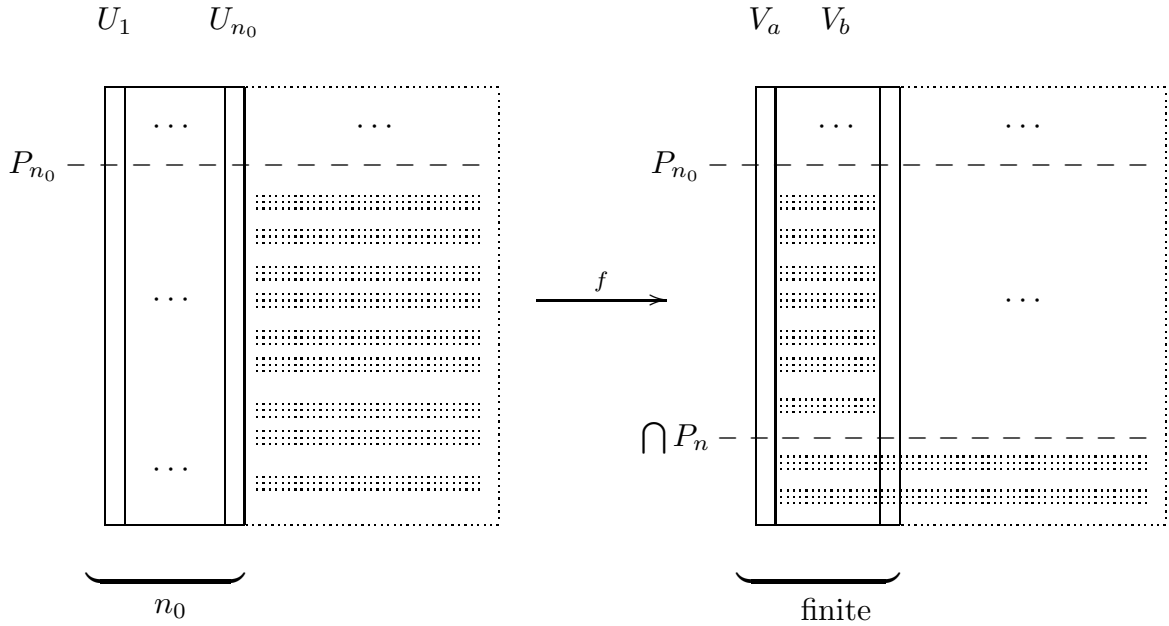
\ignore{
$$\xymatrixcolsep{0.5pc}\xymatrixrowsep{1pc}
\xymatrix{
 &U_1 && U_{n_0} &&&&&& &&&& 
 &V_a & V_b &&&\\
 &\save[0,0];[6,0]**\frm{-}\restore \save[0,0];[6,2]**\frm{-}\restore
\save[0,0];[6,8]**\frm{.}\restore &\save+<0ex,-2.5ex> \drop{\cdots}\restore 
& \save[0,0];[6,0]**\frm{-}\restore &&&\save+<0ex,-2.5ex>
\drop{\cdots}\restore &&& &&&& 
 &\save[0,0];[6,0]**\frm{-}\restore \save[0,0];[6,2]**\frm{-}\restore
\save[0,0];[6,8]**\frm{.}\restore &\save+<0ex,-2.5ex> \drop{\cdots}\restore 
& \save[0,0];[6,0]**\frm{-}\restore &&&\save+<0ex,-2.5ex>
\drop{\cdots}\restore\\ 
P_{n_0} \ar@{--}[0,9] &&&&& &&&& &&&& 
P_{n_0} \ar@{--}[0,9] &&&&& &&&&\\
 &&&&\save+<-3ex,0ex> \ar@3{.}[rrrrr] \ar@3{.}[rrrrr]<2.5ex> \restore
&&&&& &&&& 
 &\ar@3{.}[rr] \ar@3{.}[rr]<2.5ex> &&\\
 &&\save+<0ex,0ex> \drop{\cdots}\restore &&\save+<-3ex,0ex> \ar@3{.}[rrrrr]
\ar@3{.}[rrrrr]<2ex> \restore &&&&& &\ar[rrr]^f &&& 
 &\ar@3{.}[rr] \ar@3{.}[rr]<2ex> &&&&&\cdots\\
 &&&&\save+<-3ex,0ex> \ar@3{.}[rrrrr] \ar@3{.}[rrrrr]<2ex> \restore
&&&&& &&&& 
 &\ar@3{.}[rr] \ar@3{.}[rr]<2ex> \ar@3{.}[rr]<-3ex> &&\\
 &&&&\save+<-3ex,0ex> \ar@3{.}[rrrrr] \ar@3{.}[rrrrr]<2ex> \restore
&&&&& &&&& 
\bigcap P_n \ar@{--}[0,9] &&&&& &&&&\\
 &&\save+<0ex,3ex> \drop{\cdots}\restore &&\save+<-3ex,0ex>
\ar@3{.}[rrrrr]<2ex> \restore &&&&& &&&& 
 &\ar@3{.}[0,8]<1ex> \ar@3{.}[0,8]<3.5ex> &&&& &&&&\\
 &\save[0,0];[0,2]**\frm{_)}\restore &\save+<0ex,-3ex> \drop{n_0}\restore
&&&&&&& &&&& 
 &\save[0,0];[0,2]**\frm{_)}\restore &\save+<0ex,-3ex>
\drop{\txt{finite}}\restore &&
 }$$
}

Since the
summand $\bigcap_{n\in\NN} P_n\bigl( \bigoplus_J
V_j \bigr)$ on the far right of the pivotal inclusion is
a correction term which can be suppressed in most applications, the
lemma says that, after a bit of trimming on both sides, a cofinite subproduct of
the $U_i$'s maps to a finite subsum of the $V_j$'s.  

Before we give some of the arguments to illustrate the techniques, we derive
the following consequence of  Theorem 10.

\proclaim{Corollary 12} {\rm (Chase \cite{\Chase})}

 $\bullet$ If there exists a cardinal number $\aleph$ such that all of
the left $R$-modules
$R^I$ are direct sums of $\aleph$-generated submodules (this being, e.g., the
case if the projectives in $\Rmod$ are closed under direct products), then
$R$ is left perfect.   

$\bullet$ If there is a cardinal number $\aleph$ such that all left
R-modules are direct sums of $\aleph$-generated modules, then $R$ is left
artinian.
\endproclaim 

\demo{Proof of Corollary 12}  From the hypothesis of the first
assertion we deduce that $R$ has the descending chain condition for finitely
generated right ideals, those being among the matrix subgroups of the right
$R$-module $R$ by Section 3. 

Now suppose that we have the global decomposition property of the second
assertion.  Then $R$ is left perfect by the first part.  Moreover, all left
$R$-modules have the descending chain condition for annihilators of subsets of
$R$, and hence
$R$ has the ascending chain condition for annihilators of subsets of arbitrary
left $R$-modules; in other words, $R$ is left noetherian. This implies
that $R$ is indeed left artinian. \qed\enddemo

We will sketch a proof of Chase's Lemma, as the argument is clarified by
the use of
$p$-functors.

\demo{Proof of Lemma 11} The natural projection $\bigoplus_{k \in J} V_k
\rightarrow V_j$ will be denoted by $q_j$.  Assume the conclusion to be
false.  Then a standard induction yields a sequence $(n_k)_{k \in \NN}$ of
natural numbers with $n_{k+1} > n_k$, together with sequences of of pairwise
different elements $j_k \in J$, resp\. $x_k \in  P_{n_k} \bigl(\prod_{i \ge
n_k}U_i
\bigr)$, such that
$$q_{j_k}f(x_k) \notin P_{n_{k+1}}V_{j_k} \qquad\quad \text{and}
\qquad\quad q_{j_k}f(x_l) = 0 \quad \text{for}\ \ l < k.$$
Note that the definition $x = \sum_{k \in \NN} x_k \in \prod_{i \in \NN} U_i$
makes sense (indeed, in view of $x_k \in \prod_{i \ge n_k} U_i$, the sum of
the $x_k$ reduces to a finite sum in each $U_i$-component), and that for all $k
\in
\NN$ we have 
$$q_{j_k}f(x) = q_{j_k}f(x_k) + q_{j_k}f \bigl(\sum_{l >k}x_l \bigr) \ne 0,$$
because the first summand does not lie in  $P_{n_{k+1}}V_{j_k}$, whereas the
second does.  But this contradicts the fact that $f(x)$ belongs to a finite
subsum of $V_j$'s.
\qed
\enddemo

Note that the implication `$(3) \implies (2)$' of Theorem 10 is an immediate
consequence of Proposition 9.  For the converse, we refer the reader to
\cite{\ringswhose}.  We include a proof for the equivalence of $(1)$
and $(3)$ however. 

\demo{Proof of `$(1) \iff (3)$' of Theorem 10}  First suppose that $(1)$ is
satisfied, and let $P_1\supseteq P_2\supseteq
P_3\supseteq \cdots$ be a chain of $p$-functors.  Abbreviate
the intersection $\bigcap_{n\in\NN} P_n$ by $P$.  We want to prove the
existence of an index
$n_0$ such that $P_{n_0} M = PM$.  For that purpose, it is
clearly harmless to assume that $PM = 0$.  This implies that $PV = 0$ for
arbitrary direct summands $V$ of direct powers of $M$, since $P$ is in turn a
$p$-functor and, as such, commutes with direct sums and direct products.   

Set $\tau = \max\{\aleph_0, \aleph, |R|\}$, and choose a set $I$ of cardinality
at least $\tau$. By hypothesis,
$M^I \cong  \prod_{n\in\NN} M^I$ has a direct decomposition into
$\aleph$-generated summands $V_j$.  Applying Chase's Lemma to
$$\prod_{n\in\NN} M^I @>{\hphantom{xx}\text{id}\hphantom{xx}}>>
\bigoplus_{j\in J} V_j,$$ 
we obtain a natural number $n_0$ with
$$\prod_{n\ge n_0} \bigl( P_{n_0} M \bigr)^I \hookrightarrow
\bigoplus_{\text{finite}} V_j \ +\ P\bigl( \bigoplus_J V_j \bigr),$$
the final summand on the right being zero by our assumption on $M$.
But this shows the right-hand side to have cardinality at most $\tau$, while
the left-hand side has cardinality $>\tau$ if $P_{n_0}M \ne 0$.  To make the
left-hand side small enough to fit into the right, we thus need to have
$P_{n_0}M = 0$, which shows that our chain of $p$-functors becomes
stationary on $M$. Now apply Theorem 6 to obtain (3). 
	
For the converse, assume $M$ to be $\Sigma$-al\-ge\-bra\-i\-cal\-ly compact, and
set
$\aleph = \text{max}(\aleph_0, |R|)$.  It is straightforward to see that each
$\aleph$-generated submodule $U$ of $M$ can be embedded into a pure
$\aleph$-generated submodule (just close $U$ under solutions of finite
linear systems with right-hand sides in $U$, repeat the process $\aleph_0$
times, and take the union of the successive closures), the latter being a
direct summand of
$M$ by Corollary 8.  Another application of this corollary thus shows that any
maximal family of independent
$\aleph$-generated submodules, summing up to a pure submodule of $M$, must
generate all of $M$.
\qed \enddemo

As a matter of course, Theorem 10 leads us to the following answer to the
questions concerning global decompositions of modules posed in Section 2.

\proclaim{Theorem 13} {\rm (Gruson-Jensen, Huisgen-Zimmermann, [loc\.
cit.], Zimmermann \cite{\Zimakad})}
For a ring $R$, the following statements are equivalent:

{\rm (1)} Every left $R$-module is a direct sum of finitely generated
submodules. 

{\rm ($1'$)} There exists a cardinal number $\aleph$ such that every left
$R$-module is a direct sum of $\aleph$-generated submodules. 

{\rm (2)} Every left $R$-module is a direct sum of indecomposable submodules. 

{\rm (3)} Every left $R$-module is algebraically compact.
\endproclaim

\demo{Proof}  Clearly `$(1) \implies (1') \implies (3)$' by Theorem 10. The
implication `$(3)
\implies (2)$' follows from Proposition 9.

To derive `$(2) \implies (3)$' from Theorem 10, keep in mind that each left
$R$-module $M$ can be embedded as a pure submodule into an algebraically
compact module $N$.  By $(2)$ all direct products $N^I$ are direct sums of
indecomposable modules, all of which have local endomorphism rings by
Proposition 9.  Thus $N$ is $\Sigma$-al\-ge\-bra\-i\-cal\-ly compact by Theorem 10, and
so is $M$ by Corollary 8.  

`$(3) \implies (1)$'. By $(3)$, all pure inclusions of left $R$-modules split,
and therefore all left $R$-modules are pure projective.  As we mentioned
before, due to Warfield \cite{\Warf}, this means that every module $M$ is a
direct summand of a direct sum of finitely presented modules $U_i$, say $M
\oplus M' = \bigoplus_{i \in I} U_i$.  On the other hand, we already know that
$(3)$ implies decomposability of all left $R$-modules into submodules with
local endomorphism rings, and hence the Krull-Remak-Schmidt-Azumaya Theorem,
applied to $M \oplus M'$, yields the required decomposition property for $M$.
\qed
\enddemo 

In particular, we retrieve the following result of Fuller \cite{\Ful}:
Condition (1) of Theorem 13 is equivalent to

(4) Every left $R$-module $M$ has a direct sum decomposition $M =
\bigoplus_{i \in I} M_i$ which complements direct summands (i.e., given any
direct summand $N$ of $M$, there exists a subset $I' \subseteq I$ such that $M
= N \oplus \bigoplus_{i \in I'} M_i$).  

Due to the fact that every module can be purely embedded into an
algebraically compact (= pure injective) one  --  go back to the first of the
examples in Section 1  --  one can build a homology theory based on pure
injective resolutions, in analogy with the traditional homology theories. 
(Alternatively, one can use the fact that every module is an epimorphic image
of a pure projective module under an epimorphism with pure kernel and consider
pure projective resolutions.) In particular, it makes sense to speak of the
pure injective and pure projective dimensions of a module, and to define the
left pure global dimension of a ring
$R$ to be the supremum of the pure injective dimensions of its left modules. 
By playing off the two arguments of the $\Hom$-functor against each other as in
the case of the traditional global dimensions, one observes that the left pure
global dimension of a ring equals the supremum of the pure projective
dimensions of its left modules.  

Thus the rings pushed into the limelight by Theorem 10 are precisely the
ones for which the left pure global dimension is zero. Since they can be
equivalently described by the requirement that all pure inclusions of left
modules split, they are also referred to as the {\it left pure semisimple\/}
rings.  We already know from Corollary 12  that they are necessarily left
artinian, but this of course does not tell much of the story.

\head{5. Rings of vanishing left pure global dimension}\endhead

On the negative side, the rings of the title are still not completely
understood.  Ironically, however, this fact also has an upside:  Namely, as a
host of inconclusive arguments on this theme appeared in circulation, numerous
interesting insights resulted.  The purpose of this section is to describe a
representative selection of such insights and to delineate the {\it status quo}
for further work on the subject. 

 The first milestone along the way was the recognition that vanishing
of the pure global dimension on both sides takes us to a class of thoroughly
studied rings.  In fact, the solution to our problems in this left-right
symmetric situation closely parallels the outcome in the commutative case.

\proclaim{Theorem 14} {\rm (Auslander \cite{\Aus}, Ringel-Tachikawa
\cite{\RinTac}, Fuller-Reiten \cite{\FuRei})}
A ring $R$ has finite representation type if and only if the left {\rm and}
right pure global dimensions of
$R$ are zero. 
\endproclaim

The implication `only if' was shown independently by Auslander and
Ringel-Tachikawa in \cite{\Aus} and \cite{\FuRei}; the following easy
argument is due to Zimmermann [unpublished].  The converse was first
established by Fuller-Reiten; we will re-obtain it as a consequence of the
versatile Tool 17 at the end of this section.

\demo{Proof of `only if'}  Suppose that $R$ has finite representation type,
and set $M = \bigoplus_{1 \le i \le n}M_i$, where $M_1, \dots, M_n$
represent the isomorphism types of the indecomposable finitely generated
left $R$-modules. The ring $R$ being left artinian, the module $M$ has finite
length, whence its endomorphism ring is semiprimary and so, in particular,
left perfect.  We denote the opposite of this endomorphism ring by $S$, thus
turning $M$ into a right $S$-module.  To prove that the left pure global
dimension of $R$ is zero, it   suffices to check that every left
$R$-module
$A$ is a direct sum of copies of the
$M_i$.  For that purpose, we start by choosing a pure projective presentation of
$A$.  Since in our present situation the pure projective modules are precisely
the objects in
$\text{Add}(M)$ (see the remark following the definition of algebraic
compactness in Section 1), this amounts to the existence of a short exact
sequence 
 $$ \qquad 0 \rightarrow K \rightarrow M^{(I)} \rightarrow A
\rightarrow 0, \tag \dagger$$
such that $K$ is pure in $M^{(I)}$.  Our goal is to show that this
sequence splits.  The module $M^{(I)}$ being a direct sum of finitely generated
modules with local endomorphism rings, the Crawley-J\'onsson-Warfield
theorem (see, e.g.,\cite{\AndFul}) will then
tell us that
$A$ is in turn a direct sum of copies of the $M_i$.

To show splitness of $(\dagger)$, we observe that the pure projectivity of $M$
guarantees the following sequence of left $S$-modules to be exact: 
$$ \quad  0 \rightarrow \Hom_R(M,K) \rightarrow \Hom_R(M,M^I)
\rightarrow \Hom_R(M,A) \rightarrow 0 \tag \ddagger$$
  In fact, this sequence is even pure
exact.  We will deduce this from the elementary fact that the pure exact sequences
are precisely those short exact sequences whose exactness is preserved by
all functors
$\Hom_S(B,-)$ with finitely presented first argument $B$.  So let $B$ be a
finitely presented left $S$-module, and consider the following commutative
diagram, the columns of which reflect the adjointness of Hom and tensor
product; for compactness,
$\Hom$-groups are denoted by square brackets:

\ignore{
$$\xymatrixcolsep{3.2pc}\xymatrixrowsep{2pc}
\xymatrix{
0 \ar[r] &[B,[M,K]] \ar[r] \ar[d]_{\cong} &[B,[M,M^{(I)}]] \ar[r]
\ar[d]_{\cong} &[B,[M,A]] \ar[r] \ar[d]^{\cong} &0\\
0 \ar[r] &[M\otimes_S B,\, K] \ar[r] &[M\otimes_S B,\, M^{(I)}] \ar[r]
&[M\otimes_S B,\, A] \ar[r] &0
}$$
}

\noindent The lower row is exact, since $M \otimes_S B$ is a finitely presented
$R$-module, and pure projective as such.  Consequently, the upper row is
exact as well.

The sequence $(\ddagger)$ thus provides us with a pure inclusion of left
$S$-modules,
$$0 \rightarrow \Hom_R(M,K) \rightarrow \Hom_R(M,M^{(I)}) \cong S^{(I)}.$$
Since $S$ is left perfect, this sequence actually splits by Observations 0 of
Section 1.  The splitness of
$(\dagger)$ can now be gleaned from the following commutative diagram which
results from the fact that
$M$ is a generator for $\Rmod$ and hence provides us with a functorial isomorphism
$M \otimes_S \Hom_R(M,-) \rightarrow \text{id}_{\Rmod}$:

\ignore{
$$\xymatrixcolsep{3.2pc}\xymatrixrowsep{2pc}
\xymatrix{
0 \ar[r] &M\otimes_S [M,K] \ar[r] \ar[d]_{\cong} &M\otimes_S [M, M^{(I)}]
\ar[r] \ar[d]_{\cong} &M\otimes_S [M,A] \ar[r] \ar[d]^{\cong} &0\\
0 \ar[r] &{\quad K\quad} \ar[r] &{\quad M^{(I)}\quad} \ar[r] &{\quad A\quad}
\ar[r] &0 
}$$
} 

\noindent Indeed, splitness of the upper row yields splitness of the lower. \qed
\enddemo
\bigskip

\proclaim{Still open: The pure semisimplicity problem} Is every ring with
vanishing one-sided pure global dimension of finite representation type?
\endproclaim

In view of the preceding theorem, the question can be rephrased as to whether
pure semisimplicity is a left-right symmetric property. 
In the sequel, we list a number of partial results which resolve the
problem in the positive for classes of rings exhibiting some  --  even if faint 
--  `commutativity symptoms'.  Subsequently, we will see that in general the
left pure semisimple rings at least come very close to having finite
representation type, in a sense to be made precise.

The first statement of the next theorem is due to
Auslander
\cite{\Auslarge}, while the second strengthened version was proved by Herzog
\cite{\Herz}.  The third assertion was first obtained by Herzog and
later derived by Schmidmeier from a more general duality principle (see
\cite{\Herz} and \cite{\Schmid}).  

\proclaim{Theorem 15} The answer to the pure semisimplicity question is `yes'
within the following classes of rings:

$\bullet$ Artin algebras.

$\bullet$ More generally, rings with self-duality.

$\bullet$ P.I. rings.    
\endproclaim

We will sketch an elementary proof for a slight generalization of the first
assertion.  It relates left-right symmetry of pure semisimplicity directly to
the existence of almost split maps.  Namely, as the author showed in 
\cite{\Huistrong}, the following is true:\smallskip 

$\bullet$ If $\rmod$ has left almost split maps, then vanishing of
the left pure global dimension of $R$ implies finite representation type.

\demo{Proof of the preceding statement}  The following argument is
inspired by  \cite{\Ausmal}, where the notions of a preprojective/preinjective
partition are introduced.  Suppose that
$R$ has left pure global dimension zero.  By Corollary 12, this forces $R$ to
be left artinian.  Moreover, suppose that each indecomposable object
$A$ in $\rmod$ is the source of a left almost split map, i.e., of a
nonsplit monomorphism $\phi: A \rightarrow B$ such that each homomorphism from
$A$ to another object in $\rmod$, which is not a split monomorphism, factors
through
$\phi$. Our proof for representation-finiteness of $R$ uses Tool 17 below, as
well as the following well-known duality discovered by
Hullinger
\cite{\Hull} and Simson \cite{\Sim}:  If $E$ is the minimal injective
cogenerator for $\Rmod$ and $T$ the opposite of the endomorphism ring of $E$,
then $T$ is twosided artinian, again has left pure global dimension zero, and
the functor $\Hom_R(-,E)$ induces a Morita duality $\rmod \rightarrow
\text{mod-}T$.  As in the proof of Corollary 18, we exploit the descending chain
condition for finite matrix subgroups, satisfied globally in $T\text{-mod}$, to
obtain the ascending chain condition for finite matrix subgroups in arbitrary
right $T$-modules.  The remark at the end of Section 3 now shows
that each direct sum of finitely generated right $T$-modules is in fact
noetherian over its endomorphism ring. 

Let
$\DD$ be a transversal of the finitely generated indecomposable right
$T$-modules.  In a first step we will establish a {\it strong
preprojective partition} on $\DD$, namely an
ordinal-indexed partition $\DD = \bigcup_{\alpha} \DD_{\alpha}$ of
$\DD$ into pairwise disjoint finite subsets $\DD_{\alpha}$ such that each
$\DD_{\alpha}$ is a minimal generating set for 
$\DD \setminus \bigcup_{\beta < \alpha} \DD_{\beta}$ and such that, for each
$\alpha$,
$\DD_{\alpha}$  consists precisely of those objects in
$\DD \setminus \bigcup_{\beta < \alpha} \DD_{\beta}$ for which every
epimorphism
$\bigoplus_{finite} D_i \rightarrow D$ with $D_i \in \DD \setminus
\bigcup_{\beta < \alpha} \DD_{\beta}$ splits.  This claim will follow from an
obvious transfinite induction if we can show that every nonempty subset
$\DD' \subseteq \DD$ contains a finite generating set, and that any minimal such
generating set
$\DD_0$ consists precisely of those $D \in \DD'$ which have the property
that all epimorphisms from
$\add \DD'$ onto $D$ split.  For simplicity, we denote the subset $\DD'$ again
by $\DD$.  Let $M = \bigoplus_{i \in I} M_i$ be the direct sum of the
objects in $\DD$, and $S$ the endomorphism ring of $M$.  Since $M$ is
noetherian over $S$, we obtain a finite subset $I' \subseteq I$ such that $M =
S\bigl(\bigoplus_{i \in I'} M_i\bigr)$; this shows that the set $\{M_i \mid i
\in I'\}$ generates $\DD$.  Let $\DD_0 \subseteq \DD$ be a minimal finite
generating set for $\DD$.  Then, clearly, each object $D \in \DD$, with the
property that arbitrary epimorphisms from $\add \DD$ onto $D$ split, belongs
to $\DD_0$.  For the converse, let $D \in \DD_0$ and 
$f: X \rightarrow D$ be an epimorphism with $X \in \add\DD$.  Moreover choose an
epimorphism
$g:D^l \oplus \bigoplus_{1 \le i \le m}D_i \rightarrow X$, where the $D_i$ are
objects in $\DD_0\setminus \{D\}$. Note that splitting of $h = fg$ will
imply splitting of $f$.  In order to see that $h$ splits, denote the
endomorphism ring of $D$ by $S(D)$, and let $h_1, \dots ,h_l$ be the
restrictions of $h$ to the various copies of $D$ occurring as summands of the
domain of $h$.  Assume
that all of the $h_i$ are non-isomorphisms.  The ring $S(D)$ being local, this
means that the $h_i$ belong to the radical $J(D)$ of $S(D)$.  Let
$D' \subseteq D$ be the trace of $\bigoplus_{1 \le i \le m}D_i$ in $D$. The
fact that $J(D)D$ is superfluous in $D$ (keep in mind that $D$ is noetherian
over $S(D)$), applied to the equality
$D = h_1 D + \dots + h_l D + D'$ yields $D = D'$.  But this means that
$\DD_0 \setminus \{D\}$ is still a generating set for $\DD$,
a contradiction to our minimal choice of $\DD_0$.  Hence one of the $h_i$ is an
isomorphism, and consequently $h$ splits.

Next we apply the duality $\Hom_R(-,E)$ to $\DD$, to obtain a {\it
strong preinjective partition} of the transversal $\CC = \Hom_R(\DD,E)$ of the
indecomposable objects in $\rmod$.  In other words, we obtain a partition $\CC =
\bigcup_{\alpha} \CC_{\alpha}$, such that each $\CC_{\alpha}$
is a minimal finite cogenerating set for $\CC \setminus \bigcup_{\beta <
\alpha} \CC_{\beta}$, where the $\CC_{\alpha}$ are pinned down by the following
additional property $(\dagger)$:  Namely, each $\CC_{\alpha}$ consists precisely
of those objects $C$ in $\CC \setminus \bigcup_{\beta < \alpha} \CC_{\beta}$ for
which all monomorphisms from $C$ to $\add\bigl(\CC \setminus
\bigcup_{\beta <
\alpha} \CC_{\beta}\bigr)$ split.

The smallest ordinal $\tau$ such that $\CC_{\tau} = \varnothing$ is
called the length of the partition; clearly $\CC_{\beta} =
\varnothing$ for all $\beta > \tau$.  Observe that $\tau$ is a successor
ordinal.  Indeed, if $S_1, \dots, S_n$ are the simple modules in $\CC$  --
say $S_i \in \CC_{\beta_i}$  --  then $\beta = 1 + \max(\beta_1,
\dots, \beta_n)$ is the length of the partition.  This is due to the fact that
each object $C \in \CC_{\beta}$ is non-simple indecomposable and hence gives
rise to a nonsplit embedding of one of the $S_i$ into $C$; but such a nonsplit
monomorphism is incompatible with property $(\dagger)$ of $\CC_{\beta_i}$.

Finally, we show that $\tau$ is bounded above by the first infinite
ordinal number $\omega$.  In light of the preceding paragraph, this will yield
finiteness of $\tau$ and thus finiteness of $\CC$, i.e., finitenenss of the
representation type of $R$.  Suppose to the contrary that there
exists an object $A \in \CC_{\omega}$, and let $\phi: A \rightarrow B$ be a left
almost split monomorphism with $B \in \rmod$.  Moreover, for $\alpha < \omega$,
let
$B^{\alpha}$ be the reject of $\CC_{\alpha}$ in $B$, and observe that
$B^{\beta} \subseteq B^{\alpha}$ whenever $\beta < \alpha <
\omega$.  Due to the finite length of $B$, the resulting chain of rejects
becomes stationary, say at $\gamma < \omega$.  This shows that
$B/B^{\gamma}$ is cogenerated by $\CC_{\alpha}$ for all ordinal numbers
$\alpha$ between $\gamma$ and $\omega$, which places a copy of $B/B^{\gamma}$
into
$\add\bigl(\bigcup_{\beta \ge \omega} \CC_{\beta}\bigr)$.  On the other hand,
since each of the $\CC_{\alpha}$ for $\alpha < \omega$ cogenerates $A$  -- 
recall that $A \in \CC_{\omega}$  --  we obtain a family of non-split
monomorphisms $f_{\alpha}: A \rightarrow C_{\alpha}$, where $C_{\alpha} \in
\add\CC_{\alpha}$ and
$\alpha$ again runs through the ordinals between $\gamma$ and $\omega$.  Now
left almost splitness of  $\phi$ permits us to factor all of these maps through
$\phi$, say
$f_{\alpha} = g_{\alpha} \phi$ for a suitable homomorphism $g_{\alpha}: B
\rightarrow C_{\alpha}$.  By the choice of $\gamma$, the reject $B^{\gamma}$ is
contained in the kernel of each $g_{\alpha}$, which shows that the homomorphism
$\overline
\phi: A \rightarrow B/B^{\gamma}$ obtained by composing $\phi$ with the
canonical map $B \rightarrow B/B^{\gamma}$ is still a monomorphism.  It is in
turn nonsplit because $\phi$ does not split.  But in view of the above placement
of
$B/B^{\gamma}$ relative to our preprojective partition, this contradicts
property $(\dagger)$ of $\CC_{\omega}$.  Hence the assumption that
$\CC_{\omega}$ be nonempty is absurd and our argument is complete. \qed  
\enddemo 

In the meantime, Simson has replaced the pure semisimplicity problem by a
conjecture, to the effect that the answer is
negative (see \cite{\SimArt, \Simcounter}).  His key to a potential class of 
counterexamples is the following:

\proclaim{Connection with a strong Artin problem for division rings} The
answer to the pure semisimplicity question is positive if and only if the
following is true:

For every simple $D$-$E$-bimodule $M$, where $D$ and $E$ are division rings, such
that $\dim (_D M) < \infty$ and $\dim(M_E) = \infty$, there exists a
non-finitely generated indecomposable left module over the triangular matrix ring
$$\pmatrix D & M \\ 0 & E \endpmatrix.$$
\endproclaim

For a proof see \cite{\SimArt}.

On the other hand, the next theorem guarantees the rings of
vanishing left pure global dimension to at least have a very sparse supply of
indecomposable finitely generated modules.  The result was independently proved
by Prest 
\cite{\Prest} and Zimmermann and the author  \cite{\HuiZimII}. In particular,
it relates the pure semisimplicity conjecture to the quest for a better
understanding of those left artinian rings of infinite representation type
which fail to satisfy the conclusion of the second Brauer-Thrall Conjecture. 
Recall that for a left artinian ring
$R$ of infinite cardinality, the latter postulates the following:  If
$R$ has infinite representation type, then there exist infinitely many distinct
positive integers $d_n$ such that, for each $n$, there are infinitely many
isomorphism classes of indecomposable left
$R$-modules of composition length $d_n$.  While this conjecture has long been
confirmed for finite dimensional algebras over algebraically closed fields,
in \cite{\Rinspecies} Ringel  constructed a class of artinian rings with
infinite center which violate this implication.

\proclaim{Theorem 16} 
Suppose that $R$ has left pure global dimension zero. Then:

$\bullet$ For each $d\in\NN$, there exist only finitely many left $R$-modules
of length $d$, up to isomorphism.

$\bullet$ For each $d\in\NN$, there exist only finitely many 
length-$d$ modules in $\finpresR$, up to isomorphism.
\endproclaim

The obstacle one meets in trying to strengthen the latter assertion for
\underbar{right}
$R$-modules to the level of that for the left lies in the fact that
it is not known whether left pure semisimple rings are necessarily right
artinian.  In fact, it is known that a positive answer to this question would
resolve the pure semisimplicity problem in the positive.  

The following duality for finite matrix subgroups, proved in
\cite{\HuiZimII} and  -- in model-theoretic terms  --  also in
\cite{\Prest}, is one of the pivotal tools in proving Theorem 16. We
include it here because it yields some by-products of independent interest.
\smallskip 

\proclaim{Tool 17} If $M$ is an $R$-$S$-bimodule and $C$ an injective
cogenerator for Mod-$S$, then the following lattices are anti-isomorphic:
  
$\bullet$ the lattice of finite matrix subgroups of the left $R$-module $M$

\noindent and 

$\bullet$ the lattice of finite matrix subgroups of the right $R$-module
Hom$_S(M,C)_R$.
\endproclaim

\demo{Sketch of proof}
Keep in mind that 
the finite matrix subgroups of $_RM$ are precisely the kernels of the maps 
$M \rightarrow Q \otimes_R M$, $m \mapsto q \otimes m$ with $Q$ finitely
presented and $q \in Q$, and that they can alternately be given the form
$\Hom_R(P,M)(p)$ with $P$ finitely presented and $p \in P$ (see Observations
3(1),(2)  --  a proof can be found in
\cite{\HuiZimII, Lemma 1}).

First one checks that, given a finite matrix subgroup $U = \Hom_R(P,M)(p)$ of
$_RM$ with $P$ and $p$ as above, the $C$-dual
$\Hom_S(M/U, C)$ is a finite matrix subgroup of the right $R$-module
$M^+ = \Hom_S(M, C)$.  From the exactness of the sequence
$$\Hom_R(P,M) @>\psi>> M @>\text{canon}>> M/U \rightarrow 0,$$ 
where $\psi$ is evaluation at $p$, one obtains exactness
of the sequence of induced maps 
$$0 \rightarrow \Hom_S(M/U,C) \rightarrow \Hom_S(M,C) @>\psi^*>>
\Hom_S(\Hom_R(P,M),C).$$
Using the canonical isomorphism $\tau: \Hom_S(M,C) \otimes_R P \rightarrow
\Hom_S(\Hom_R(P,M),C)$ and the fact that the map $\tau^{-1}\psi^*$ sends $f \in
\Hom_S(M,C)$ to $f \otimes p$, one deduces the first claim.  Clearly, this
finite matrix subgroup of $M^+$ coincides with the annihilator
$\Ann_{M^+}(U)$.

Similarly, one verifies that, for any finite matrix subgroup $V$ of the right
$R$-module $M^+$, the annihilator $\Ann_M(V)$ is a finite matrix subgroup of
the left $R$-module $M$.

It is now a question of routine to check that the lattice anti-homomorphisms
$U \mapsto \Ann_{M^+}(U)$ and $V \mapsto \Ann_M(V)$ are inverse to each
other. \qed \enddemo  

It was shown by Crawley-Boevey that, among the finite dimensional algebras over
an algebraically closed field, the ones of finite representation type are
characterized by the nonexistence of generic modules, i.e., of non-finitely
generated indecomposable endofinite modules.  (As the term suggests, a
module is {\it endofinite} if it has finite length over its endomorphism ring.) 
The following consequence of the above duality shows that the absence of generic
objects is offset by the richest possible supply of endofinite modules in the
representation-finite case.

\proclaim{Corollary 18}  {\rm (\cite{\Prest, \HuiZimII})} For any ring
$R$, the following statements are equivalent:

{\rm (1)} All objects in $\Rmod$ are endofinite.

{\rm (2)} Same as {\rm (1)} for $\operatorname{Mod-}R$.

{\rm (3)} The left and right pure global dimensions of $R$ vanish.
\endproclaim  

\demo{Proof} It clearly suffices to show the equivalence
of (1) and (3).  If (1) holds, Tool 17, applied to $R$ and the opposites of
the  endomorphism rings of the modules considered, yields the descending chain
condition for finite matrix subgroups in all {\it right}
$R$-modules.  That the left modules satisfy this chain condition is immediate
from our hypothesis.  By Theorem 6, this shows that all
$R$-modules, left and right, are algebraically compact, i.e\. (3) holds.

Now assume (3).  In view of the descending chain condition for matrix subgroups,
satisfied globally for left and right $R$-modules, Tool 17  provides us with
the ascending chain condition for finite matrix subgroups as well.  Let
$M$ be any left
$R$-module and $S$ its endomorphism ring. We know that all finitely generated
$S$-submodules of $M$ are matrix subgroups, so if we can show that
all matrix subgroups of $_RM$ arise from finite matrices, we are done.  But this
finiteness condition follows easily from our hypothesis as follows:  Indeed, let
$[\A, \alpha]$ be any matrix functor on $\Rmod$, where $\A = (a_{ij})_{i \in I, j
\in J}$.  Moreover, for any finite subset $I'$ of $I$, let $\A(I')$ be the
matrix consisting of those rows of $\A$ which are indexed by $I'$.  Since $\A$
is row-finite, the matrix $\A(I')$ is actually finite in effect.  Now let
$I_1$ be any finite subset of $I$ and note that the finite matrix subgroup
$[\A(I_1),\alpha]M$ of $M$ contains $[\A, \alpha]
M$.  If the inclusion is proper, there exists a finite subset $I_2 \subseteq I$
containing $I_1$ such that $[\A(I_2), \alpha]M$ is properly contained in
$[\A(I_1), \alpha]M$.  If $[\A, \alpha]M$ is still strictly contained in
$[\A(I_2), \alpha]M$, we repeat the process.  Eventually, it has to
terminate by our chain condition, which shows that $[\A, \alpha]M$ equals some
finite matrix subgroup $[\A(I_m),\alpha]M$, and our argument is complete. \qed
\enddemo

We conclude this section by looping back to its beginning and filling in the
proof of the as yet unjustified implication of Theorem 14. We prepare with the
following proposition (not needed in full strength here) which generalizes
Crawley-Boevey's observation that an endofinite direct sum of modules with
local endomorphism rings involves only finitely many isomorphism types of
summands
\cite{\CrawBoeII, Proposition 4.5}.

\proclaim{Proposition 19} Suppose that $(M_i)_{i\in I}$ is a family of left
$R$-modules with local endomorphism rings. Let $M= \bigoplus_{i\in I} M_i$,
and denote by
$S$ the endomorphism ring of $M$, by $J(S)$ the Jacobson radical of $S$.

If there exists a natural number $N$ with $J(S)^NM=0$ such that, moreover,
$J(S)^kM$ is finitely generated over $S$ for all $0\le k\le N$, then the
family $(M_i)_{i\in I}$ involves only finitely many isomorphism classes.
\endproclaim

\demo{Proof} Suppose $N$ is as in the claim, and start by noting that
$\Hom_R(M_i,M_j)\subseteq J(S)$, whenever $M_i \not\cong M_j$; this is due to
the locality of the endomorphism rings of the $M_i$. Next observe that, for
any finite sequence $x_1,\dots,x_r$ in $M$, we have
$$\sum_{i=1}^r Sx_i \subseteq \bigoplus_{i\in I'} M_i + \sum_{i=1}^r J(S)x_i,$$
where $I'\subseteq I$ is the `closure of the supports of the $x_l$, $1\le l\le
r$, under isomorphism'; by this we mean that 
$$I'= \{i\in I\mid M_i\cong M_j
\text{\ for some\ } j\in \bigcup_{1\le l\le r} \supp(x_l)\}.$$

Now pick generators $x_{k1},\dots,x_{kr_k}$ for each of the left $S$-modules
$J(S)^kM$, where $k$ runs from $0$ to $N-1$. Moreover, for each $k$, let $I_k$
be the closure of $\bigcup_{1\le l\le r_k} \supp(x_{kl})$ under isomorphism.
Then
$$\align M= \sum_{l=1}^{r_0} Sx_{0l} &\subseteq \bigoplus_{i\in I_0} M_i
+\sum_{l=1}^{r_0} J(S)x_{0l} \subseteq \bigoplus_{i\in I_0} M_i
+\sum_{l=1}^{r_1} Sx_{1l}\\
 &\subseteq \bigoplus_{i\in I_0\cup I_1} M_i +\sum_{l=1}^{r_1} J(S)x_{1l}
\subseteq \bigoplus_{i\in I_0\cup I_1} M_i +\sum_{l=1}^{r_2} Sx_{2l},
\endalign$$
and since $J(S)x_{N-1,l} =0$ for all $l$ by hypothesis, an obvious induction
yields
$$M= \bigoplus_{i\in I_0\cup \dots\cup I_{N-1}} M_i.$$
But the modules $M_i$ with $i\in I_0\cup \dots\cup I_{N-1}$ fall into finitely
many isomorphism classes by construction, and our argument is complete.
\qed\enddemo

On the side, we point out that, in Proposition 19, neither the hypothesis that
$J(S)^NM=0$ for some $N$, nor the
condition that $J(S)^kM$ be finitely generated over $S$ for all $k$, suffices
to yield the conclusion: If $R=\ZZ$ and $M= \bigoplus_{p \text{\,prime}}
\ZZ/(p)$, then clearly $J(S)=0$. An example to justify the second remark is as
follows: Let
$R$ be the Kronecker algebra and
$M_n$ the preprojective string module of dimension $2n+1$. If we set $M=
\bigoplus_{n\in
\NN} M_n$, then $J(S)^kM= \bigoplus_{n\ge k+1} M_n$ is finitely generated over
$S$ for each $k$.
\smallskip

\demo{Proof of the implication `if' in Theorem 14}  Our hypothesis being
condition (3) of Corollary 18, we see that all left $R$-modules are
endofinite. Since we already know $R$ to be artinian (Corollary 12), we
only need to show that any transversal $(M_i)_{i \in I}$ of the
indecomposable objects in $\rmod$ is finite. But this is an immediate
consequence of Proposition 19. \qed
\enddemo     

For the reader interested in picking up the threads that are still dangling, we
list a number of additional references addressing the pure semisimplicity
problem: \cite{\Sim}, \cite{\SimI}, \cite{\Hull}, \cite{\Simcox},
\cite{\JenLenI}, \cite{\Zayed}, \cite{\Azufac}, 
\cite{\Fac}, \cite{\JenLen}, \cite{\JenZim}, \cite{\Huistrong}, \cite{\Azu}, 
\cite{\Her}, \cite{\Lara}, \cite{\Simpad}.

\head{6. Vanishing pure global dimension and product-completeness}\endhead

From Theorem 10 we know that the $\Sigma$-al\-ge\-bra\-i\-cal\-ly compact
modules $M$ are precisely those with the property that all direct products
of copies of
$M$ are direct sums of indecomposable components with local endomorphism
rings.  This naturally raises the question as to the structure of the
indecomposable summands in such decompositions.

More generally, we ask:  {\it What
are the indecomposable summands of large 
$\Sigma$-al\-ge\-bra\-i\-cal\-ly compact direct products $\prod_{i\in I} M_i$?}

This problem is
further motivated by the following points: (a) The question is intimately
related to the pure semisimplicity problem.  The connection first surfaced in a
theorem of Auslander (Theorem 20 below), and is reinforced by results of Gruson,
Garavaglia, and Krause-Saor\'\i n (also compare with Corollary 18). And (b),
direct powers of non-generic $\Sigma$-al\-ge\-bra\-i\-cal\-ly compact modules
provide prime hunting ground for generic objects, as evidenced by Theorems
2.11 and 2.12 (due to Krause and Ringel, respectively) in Zwara's contribution
to this volume.  

Slightly extending the terminology of Krause-Saor\'\i n \cite{\KraSao},
we call a family $(M_i)_{i \in I}$ {\it product-complete}, if the direct
product $\prod_{i \in I} M_i$ belongs to $\Add(\bigoplus_{i \in I}
M_i)$, the closure of $\{M_i \mid i \in I\}$ in $\Rmod$ under formation of
direct summands and direct sums.  As in \cite{\KraSao}, a module $M$ will be
labeled product-complete if all families consisting of copies of $M$ have this
property.

According to Gruson \cite{\Gru} and Garavaglia
\cite{\Gara}, an indecomposable module $M$ is product-complete if and only if
it is endofinite.  This equivalence
was generalized by Krause-Saor\'\i n as follows:  For an arbitrary module $M$,
finite endolength is equivalent to the postulate that all direct summands of
$M$ be product-complete \cite{\KraSao, Theorem 4.1}.  

This yields the following hierarchy within the class of algebraically compact
modules:  

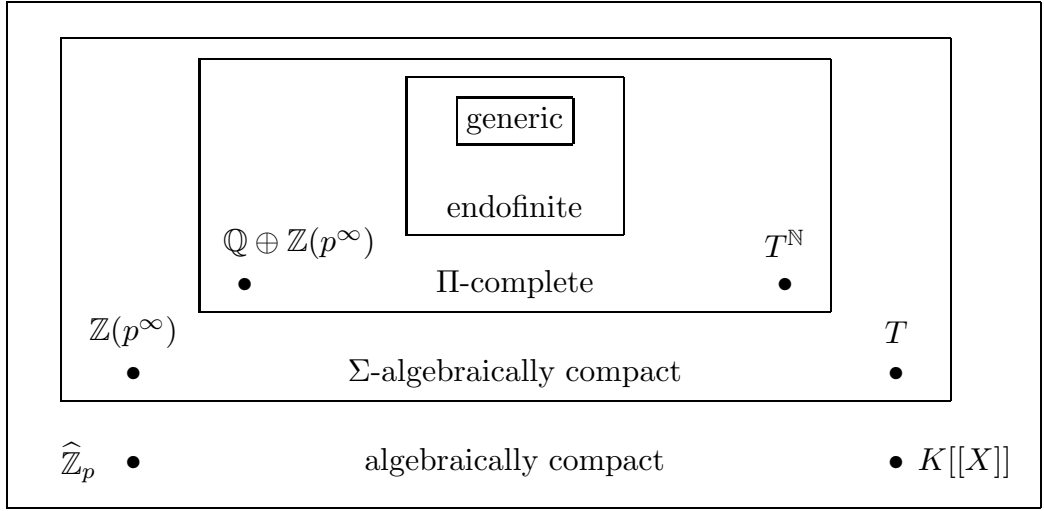
\ignore{
$$\xymatrixcolsep{2pc}\xymatrixrowsep{1.1pc}
\xymatrix{
\save[0,0]+(-8,0);[4,4]+(6,-3)**\frm{-}\restore 
\save[0,0]+(-14,4);[5,4]+(16,-5)**\frm{-}\restore &&&&\\
 & \save[0,0]+(-5,6);[2,2]+(5,-3)**\frm{-}\restore &\txt{generic}
\save**\frm{-}\restore \save[0,0]+(-12,4);[1,0]+(12,-4)**\frm{-}\restore &&\\
 &&\txt{endofinite} &&\\
 &\bullet \save+<4ex,3ex> \drop{{\Bbb Q}\oplus{\Bbb Z}(p^{\infty})}\restore
&\Pi\hbox{-complete} &\bullet \save+<0ex,3ex> \drop{T^{\Bbb N}}\restore &\\
\bullet \save+<0ex,3ex> \drop{{\Bbb Z}(p^{\infty})}\restore
 &&\Sigma\hbox{-al\-ge\-bra\-i\-cal\-ly compact} &&\bullet \save+<0ex,3ex>
\drop{T}\restore\\
\bullet \save+<-4ex,0ex> \drop{\widehat{\Bbb Z}_p}\restore
&&\hbox{algebraically compact} &&\bullet \save+<5ex,0ex>
\drop{K[[X]]}\restore  
}$$
}
\smallskip

\noindent Each of the inclusions is proper: The algebra of power series $K[[X]]$
is algebraically compact without being $\Sigma$-al\-ge\-bra\-i\-cal\-ly compact;
and the same is true for the ring $\widehat{\ZZ}_p$ of $p$-adic integers,
viewed either as  $\ZZ$- or  $\widehat{\ZZ}_p$-module. The ring
$T = K[X_i \mid i \in \NN]/(X_i \mid i \in \NN)^m$, where
$K$ is a field and $m$ some integer $\ge 2$, is
$\Sigma$-al\-ge\-bra\-i\-cal\-ly compact (c.f\. examples following Theorem
6), but not product-complete; indeed, $T^{\NN}$ is not projective, as $T$ fails
to be coherent.  On the other hand, Schulz proved the non-projective
summands in arbitrary direct products of copies of $T$ to be all isomorphic to
the unique simple $T$-module (see \cite{\Schulz}), whence $T^{\NN}$ is
product-complete without being endofinite.  Furthermore, the direct sum of the
group of rational numbers and the Pr\"ufer group ${\Bbb Z}(p^{\infty})$ for
some prime $p$ is a product-complete abelian group, but fails to be
endofinite.  Over any Artin algebra
$R$, finally, all objects in $\rmod$ are endofinite without being
generic.  

As is to be expected in light of our previous discussion, global
product-completeness of families of $R$-modules occurs only rarely. More
precisely, we have:      

\proclaim{Theorem 20} {\rm (Auslander \cite{\Auslarge})} For any ring $R$, the
following statements are equivalent:

{\rm (1)} $R$ has finite representation type.

{\rm (2)} All
families of finitely generated indecomposable left $R$-modules are
product-com\-plete. \endproclaim

We give a proof relying only on the results established in the previous
sections.  Note that the first part of our argument
is very similar to the reasoning of Krause-Saor\'\i n (\cite{\KraSao,
3.8}).   

\demo{Proof}  `$(1) \implies (2)$'.  Given $(1)$, we know from Theorem 14 that
the left pure global dimension of $R$ is zero.  Let $(M_i)_{i \in I}$ be a
family of indecomposable objects in $\rmod$.  To see that each indecomposable
direct summand of $\prod_{i \in I} M_i$ is isomorphic to one of the $M_j$, it
suffices to show this for the case where all $M_i$ are pairwise
isomorphic since, by hypothesis, our family contains only finitely many
isomorphism types.  So suppose $M_i \cong M$ for all $i$.  Denote the
endomorphism ring of $M$ by $S$.  From Corollary 18, we know that $M$ has finite
length over $S$, and hence
$S$ is left artinian:  Indeed, letting $m_1, \dots, m_n$ be a generating set
for $M$ over
$R$, we see that $S$ embeds into $M^n$ as a left $S$-module, via $s
\mapsto (sm_k)$.  In particular, this implies that $S^I$ is projective as a
right $S$-module.  The ring $S$ being local, this means that $S^I$ is free, say
$S^I \cong S^{(J)}$.  Using the fact that the tensor functor $- \otimes_S M$
commutes with direct products in our setting, we see that $M^I \cong S^I
\otimes_S M \cong S^{(J)} \otimes M \cong M^{(J)}$, which yields $(2)$. 

For the converse, assume that $(2)$ is true, and let $(M_i)_{i \in I}$
be a transversal of the isomorphism types of the indecomposable objects in
$\rmod$.  We will obtain finiteness of $I$ by showing that the direct
product $\prod_{i \in I} M_i$ equals the direct sum $\bigoplus_{i \in I}
M_i$.  We start by observing that $\bigoplus_{i \in I} M_i$ is
algebraically compact:  Indeed, all direct powers of $\prod_{i \in I}
M_i$ are direct sums of copies of the $M_i$ by hypothesis, whence $\prod_{i \in
I} M_i$ is $\Sigma$-algebraically compact by Theorem 10; but by Theorem 6, this
is tantamount to $\Sigma$-algebraic compactness of $\bigoplus_{i
\in I} M_i$.  In particular, the pure inclusion $\bigoplus_{i
\in I} M_i \subseteq  \prod_{i \in I} M_i$ splits, say $\prod_{i \in
I} M_i = \bigoplus_{i \in I} M_i \oplus N$, and $N$ is in
turn a direct sum of copies of $M_i$'s.  If $N$ were nonzero,
we could thus find a direct summand isomorphic to some $M_k$ in $N$. 
On the other hand, we may cancel $M_k$ from the above product-sum equality to
obtain $\prod_{i \in I \setminus \{k\}} M_i \cong \bigoplus_{i \in I \setminus
\{k\}} M_i \oplus N$ (keep in mind that $\End(M_k)$ is local).  This makes
$M_k$ a direct summand of $\prod_{i
\in I \setminus \{k\}} M_i$, thus contradicting our
hypothesis. We conclude
$N=0$, which forces $I$ to be finite as required. \qed 
\enddemo

The following result, due to Krause and Saor\'\i n \cite{\KraSao,
Proposition 4.2}, characterizes the product-complete modules in terms of their
matrix subgroups.  It continues the line of Theorem 10, where we related
direct sum decompositions of direct products $M^I$ to finiteness conditions on
the lattice of matrix subgroups of $M$. 

\proclaim{Proposition 21} An object $M \in \Rmod$ is product-complete if
and only if $M$ has the descending chain condition for (finite) matrix
subgroups and all (finite) matrix subgroups of $M$ are finitely generated over
the endomorphism ring of $M$. \endproclaim

We include an elementary proof for one implication.

\demo{Proof of `only if'}  Suppose $M$ is product-complete.  Then, clearly, $M$
is
$\Sigma$-al\-ge\-bra\-i\-cal\-ly compact by Theorem 10, which is tantamount to $M$
satisfying the descending chain condition for matrix subgroups by
Theorem 6.  To see that any matrix subgroup $[\A, \alpha]M$ is finitely
generated over the endomorphism ring $S$ of $M$, recall that 
$$[\A, \alpha]M =
\Hom_R(M^I,M)(\underline m)$$
for some set $I$ and some element
$\underline m \in M^I$ (this was explained after the statement of Theorem 5). 
By hypothesis, $M^I$ is a direct summand of a suitable direct sum $M^{(J)}$,
and hence $[\A, \alpha]M = \Hom_R(M^{(J)},M)(m)$ for some $m = (m_j) \in
M^{(J')}$, where $J'$ is a finite subset of $J$.  This gives $[\A, \alpha]M =
\sum_{j \in J'} S m_j$ as required. \qed
\enddemo

In view of the final remark of Section 3A, a ring $R$ is right coherent
if and only if all finite matrix subgroups of the regular left module $R$ are
finitely generated as right ideals.  We can thus supplement Corollary 12 to
retrieve Chase's characterization of the rings whose left projective modules
are closed under direct products.

\proclaim{Corollary 22} {\rm (Chase \cite{\Chase})} $R$ has the property that
all direct products of projective left $R$-modules are again projective if and
only if
$R$ is left perfect and right coherent. \endproclaim

The following question appears of significantly lower importance than the one
with which we opened the section.  However, the fact that it is not yet
answered shows the lacunary state of our present understanding of direct sum
decompositions of large direct products.
\smallskip

\noindent {\it Given a $\Sigma$-algebraically compact module $M$, is $M^{\NN}$
product-complete?}\smallskip  

As was already observed by Krause and Saor\'\i n \cite{\KraSao}, there is
{\it some} power $M^I$ which is product-complete;  indeed, the fact that all
powers of $M$ split into summands of cardinalities bounded above by
$\max(|R|, \aleph_0)$ guarantees that, eventually, saturation with respect to the
appearance of new direct summands is reached.

\head {7. Concluding remarks on pure global dimension}\endhead

From the previous section we know that a ring $R$ has finite
representation type if and only if its left and right pure
global dimensions are zero. Moreover, by Theorem 15, the latter condition is
left-right symmetric for Artin algebras.  What can one say about the pure
global dimensions of
$R$ when they do not vanish?   

The question of how these invariants relate to other properties of the
ring and its module categories is largely open.  Some general
facts of interest are available, however, as well as some classes of
algebras where the connection is understood.  We include only a few results to
trigger interest in further investigation of the problem.
 
A very rough, but not unreasonable,
answer to the above question is this: ``That depends on the cardinality of
$R$."  As a first step in justifying this response, we present an insight due to
Gruson and Jensen  \cite{\GruJenfinal}; for an alternate approach, see
\cite{\KielSim}.

\proclaim{Theorem 23} The following implications hold for any ring $R$.  

{\rm (1)}  If $R$ is countable, but not of finite representation
type, then the left (right) pure global dimension of $R$ equals $1$.

{\rm (2)}  If the cardinality of $R$ is bounded from above
by $\aleph_t$ for some integer $t \ge 0$, the pure global dimensions of $R$ are
at most $t+1$.
 \endproclaim 

One can do far better for specialized classes of rings, however.  For example, as
was observed by Kie\l pi\'nski and Simson \cite{\KielSim}, if $R = S[X]$ is a
plynomial ring over a commutative ring $S$ of cardinality $\aleph_s$ in a set
$X$ of indeterminates which has cardinality $\aleph_t$ (with $s,t \ge 0$), then
the pure injective dimension of $R$ equals $\max(s,t) + 1$.  Our final result,
due to D\. Baer, Brune, and Lenzing
\cite{\BaeBruLen}, presents a smooth picture of pure homology for
hereditary algebras over algebraically closed base fields subject to certain
cardinality restrictions.  In that case, the pure global dimension mirrors the
representation type of
$R$ as follows:

\proclaim{Theorem 24}  Suppose that $R$ is a hereditary finite dimensional
algebra over an algebraically closed field of cardinality $\aleph_t$, where $t
\ge 2$ is an integer.  Then:

{\rm (1)} The left (right) pure global dimension of $R$ equals $2$ if and only
 if $R$ has tame, but infinite, representation type.

{\rm (2)}  The left (right) pure global dimension of $R$ equals $t+1$ if and
only if $R$ has wild representation type.
\endproclaim

Finally, we refer the reader to \cite{\Okoh}, \cite{\BaeLen} and
\cite{\BaeBruLen} for further instances in which the pure global dimension is
understood.  Motivated by the lectures of Benson at the 1998 conference in
Bielefeld, we conclude with the following

\proclaim{Problem} Given a finite group $G$ and a field $K$ of suitable
cardinality, how does the pure global dimension of the group algebra $KG$
reflect the representation type?
\endproclaim

\enddocument